\documentclass[11pt,a4paper]{article}
\usepackage[top=1.7cm,bottom=1.2cm,left=2.8cm,right=2.3cm,includefoot]{geometry}    
\usepackage{fancyhdr}
\usepackage{indentfirst}
\usepackage{amsmath}
\usepackage{amssymb}
\usepackage{graphicx}
\newtheorem{theorem}{Theorem}
\newtheorem{lemma}{Lemma}

\newtheorem{corollary}{Corollary}
\newtheorem{definition}{Definition}

\newtheorem{remark}{Remark}
\newcommand{\onetom}{1,2,\cdots,m}
\newcommand{\oneton}{1,2,\cdots,n}
\title{Achieving synchronization in arrays of coupled differential systems with time-varying couplings}
\author
{Xinlei Yi$^{1}$ \;
Wenlian Lu$^{1,2,3}$\; and \; Tianping Chen$^{1,4}$}
\date{\it\footnotesize
 1. School of Mathematical Sciences, Fudan University, Shanghai 200433, China\\ Email: yixinlei2008@yahoo.cn \\
 2. Centre for Computational Systems Biology, Fudan University, Shanghai 200433, China\\
 3. Centre for Scientific Computing, the University of Warwick, Coventry CV4
7AL, UK\\Email: wenlian@fudan.edu.cn\\
 4. School of Computer Science, Fudan University, Shanghai 200433, China\\Email: tchen@fuan.edu.cn}
\begin{document}
\maketitle \noindent \rule[1.5mm]{\textwidth}{0.1pt}
{\bf Abstract}

In this paper, we study complete synchronization of the complex dynamical networks
described by linearly coupled ordinary differential equation systems
(LCODEs). The coupling considered here is time-varying in both the network
structure and the reaction dynamics. Inspired by our previous paper \cite{LAJ}, the extended Hajnal diameter is
introduced and used to measure the synchronization in a general
differential system. Then we find that the Hajnal diameter of the linear
system induced by the time-varying coupling matrix and the largest Lyapunov
exponent of the synchronized system play the key roles in synchronization analysis
of LCODEs with the identity inner coupling matrix. As an application, we obtain a general sufficient condition
guaranteeing directed time-varying graph to reach consensus.
Example with numerical simulation is provided to show the effectiveness
the theoretical results.

\noindent\rule[1.5mm]{\textwidth}{0.1pt}

\section{Introduction}
Complex networks have widely been used in theoretical analysis of complex systems, such as Internet, World Wide Web, communication networks, social networks. A complex dynamical network is a large set of interconnected nodes, where each node possesses a (nonlinear) dynamical system and the interaction between nodes is described as diffusion. Among them, linearly coupled ordinary differential equation systems (LCODEs) are a large class of dynamical systems with continuous time and state.

The LCODEs are usually formulated as follows
\begin{eqnarray}\label{L}
\dot{x}^{i}(t)=f(x^{i}(t))+\sigma\sum_{j=1}^{m}l_{ij}Bx^{j}(t),~i=\onetom,
\end{eqnarray}
where $t\in\mathbb R^{+}=[0,+\infty)$ stands for the continuous time and $x^{i}(t)\in\mathbb R^{n}$
denotes the variable state vector of the $i$-th node, $f:\mathbb
R^{n}\rightarrow\mathbb R^{n}$ represents the node dynamic of the uncoupled system, $\sigma\in \mathbb R^{++}=(0,+\infty)$ denotes coupling strength, $l_{ij}\ge0$ with $i\neq j$ denotes the interaction between the two nodes and $l_{ii}=-\sum_{j\neq i}^{m}l_{ij}$, $B\in \mathbb R^{n,n}$ denotes the inner coupling matrix. The LCODEs model is widely used to describe the model in nature and engineering. For example, the authors study spike-burst neural activity and the transitions to a synchronized state using a model of linearly coupled bursting neurons in \cite{Dha}; the dynamics of linearly coupled Chua circuits are studied with application to image processing, and many other cases in \cite{Per}.

For decades, a large number of papers have focused on the dynamical behaviors of coupled systems  \cite{Ash}-\cite{Mil1}, especially the synchronizing characteristics. The word ``synchronization'' comes from Greek, in this paper the concept of local complete synchronization (synchronization for simplicity) is considered (see {\bf Definition 1}). For more details, we refer the readers to \cite{LAJ}.

Synchronization of coupled systems have attracted a great deal of attention \cite{AMR}-\cite{Wu}. For instances, in \cite{AMR}, the authors considered the synchronization of a network of linearly coupled
and not necessarily identical oscillators; in \cite{JQL}, the authors studied globally exponential synchronization
for linearly coupled neural networks with time-varying delay and impulsive disturbances.
Synchronization of networks with time-varying topologies was studied in \cite{Bel}-\cite{Sti}. For example, in \cite{Bel}, the authors proposed the global stability of total synchronization in networks with different topologies; in \cite{Sti}, the authors gave a result that the network will synchronize with the time-varying topology if the
time-average is achieved sufficiently fast.

Synchronization of LCODEs has also been addressed in \cite{Lu1}-\cite{Liubo1}. In \cite{Lu1}, mathematical analysis was presented on the synchronization phenomena of LCODEs with a single coupling delay; in \cite{Lu},  based on geometrical analysis of the synchronization manifold, the authors proposed a novel approach to investigate the stability of the synchronization manifold of coupled oscillators; in \cite{Liubo1}, the authors proposed new conditions on synchronization of networks of linearly coupled dynamical
systems with non-Lipschitz right-handsides. The great majority of research
activities mention above all focused on static networks whose connectivity and coupling strengths are
static. In many applications, the interaction between individuals may change dynamically. For example, communication links between agents may be unreliable due to disturbances and/or subject to communication range limitations.

In this paper, we consider synchronization of LCODEs with time-varying coupling. Similar to \cite{Lu1}-\cite{Liubo1}, time-varying coupling will be used to represent the interaction between individuals. In \cite{LAJ,LAJ1}, they showed that the Lyapunov exponents of the synchronized system and the Hajanal diameter of the variational equation play key roles in the analysis of the synchronization in the discrete-time dynamical networks. In this paper, we extend these results to the continuous-time dynamical network systems. Different from \cite{Bel1,Sti}, where synchronization of fast-switching systems was discussed, we focus on the framework of synchronization analysis with general temporal variation of network topologies. Additional contributions of this paper are that we explicitly show that (a) the largest projection Lyapunov exponent of a system is equal to the logarithm of the Hajanal diameter, and (b) the largest Lyapunov exponent of the transverse space is equal to the largest projection Lyapunov exponent under some proper conditions.

The paper is organized as follows: in Section 2, some necessary definitions, lemmas, and hypotheses are given; in Section 3, synchronization of generalized coupled differential systems is discussed; in Section 4, criteria for the synchronization of LCODEs are obtained; in Section 5, we obtain a sufficient condition ensuring directed  time-varying graph reaching consensus; in Section 6, example with numerical simulation is provided to show the effectiveness of the theoretical results; the paper is concluded in Section 7.

\noindent {\bf Notions}:
$e^{n}_{k}=[0,0,,\cdots,0,1,0,\cdots,0]^{\top}\in\mathbb R^{n}$ denotes the $n$-dimensional vector with all components zero except the $k$-th component 1, $\bf 1_{n}$ denotes the
$n$-dimensional column vector with each component 1; For a set in some Euclidean space  $U$,  $\bar{U}$ denotes the closure of $U$, $U^{c}$ denotes the complementary set of $U$, and $A\setminus B=A\cap B^{c}$; For $u=[u_{1},\cdots,u_{n}]^{\top}\in\mathbb R^{n}$, $\|u\|$
denotes some vector norm, and for any matrix $A=(a_{ij})\in\mathbb R^{n,m}$, $\|A\|$ denotes some matrix norm induced by vector norm, for example, $\|u\|_{1}=\sum_{i=1}^{n}|u_{i}|$ and $\|A\|_{1}=\max_{j}\sum_{i=1}^{n}|a_{ij}|$; For a matrix $A=(a_{ij})\in\mathbb R^{n,m}$, $|A|$ denotes a
matrix with $|A|=(|a_{ij}|)$; For a real matrix $A$, $A^{\top}$ denotes its transpose and
for a complex matrix $B$, $B^{*}$ denotes its conjugate transpose; For a set in some Euclidean space  $W$, $\mathcal O(W,\delta)=\{x:~{\rm dist}(x,W)<\delta\}$, where
${\rm dist}(x,W)=\inf_{y\in W}\|x-y\|$; $\# J$ denotes the cardinality of set $J$; $\lfloor z\rfloor$ denotes the floor function, i.e., the
largest integer not more than the real number $z$;
$\otimes$ denotes the Kronecker product.

\section{Preliminaries}

In this section we will give some necessary definitions, lemmas, and hypotheses.
Consider the following general coupled differential system:
\begin{eqnarray}
\dot{x}^{i}(t)=f^{i}(x^{1}(t),x^{2}(t),\cdots,x^{m}(t),t),~i=\onetom,\label{CDS}
\end{eqnarray}
with initial state $x(t_0)=[{x^{1}(t_0)}^{\top},\cdots,{x^{m}(t_0)}^{\top}]^{\top}\in\mathbb R^{nm}$, where $t_0\in\mathbb R^{+}$ denotes the initial time, $t\in\mathbb R^{+}$ denotes the continuous time and $x^{i}(t)=[x^{i}_{1}(t),\cdots,x^{i}_{n}(t)]\in\mathbb R^{n}$
denotes the variable state of the $i$-th node, $i=\onetom$.

For the functions $f^{i}:\mathbb
R^{nm}\times \mathbb R^{+}\rightarrow\mathbb R^{n}$, $i=\onetom$, we make
following assumption:

{\bf Assumption ${\mathbf A1}$} (a)
There exists a function $f:\mathbb R^{n}\rightarrow\mathbb
R^{n}$ such that $f^{i}(s,s,\cdots,s,t)=f(s)$ for all
$i=\onetom$, $s\in\mathbb R^{n}$, and $t\ge 0$; (b) For any $t\ge 0$, $f^{i}(\cdot,t)$ is $C^{1}$-smooth for all
$x=[{x^{1}}^{\top},\cdots,{x^{m}}^{\top}]^{\top}\in\mathbb R^{nm}$, and by $DF^{t}(x)=(\frac{\partial
f^{i}}{\partial x^{j}}(x,t))_{i,j=1}^{m}\in\mathbb R^{nm,nm}$ denotes the
Jacobian matrix of
$F(x,t)=[f^{i}(x,t)^{\top},\cdots,f^{m}(x,t)^{\top}]^{\top}$ with
respect to $x\in\mathbb R^{nm}$; (c) There exists a
locally bounded function $\phi(x)$ such that
$\|DF^{t}(x)\|\le\phi(x)$ for all $(x,t)\in\mathbb
R^{nm}\times\mathbb R^{+}$; (d) $DF^{t}(x)$ is {\it uniformly locally Lipschitz continuous}:
there exists a locally bounded function $K(x,y)$ such that
\begin{eqnarray*}
\|DF^{t}(x)-DF^{t}(y)\|\le K(x,y)\|x-y\|
\end{eqnarray*}
for all $t\ge 0$ and $x,y\in\mathbb R^{nm}$; (e) $f^{i}(x,t)$ and $DF^{t}(x)$ are both measurable for $t\ge
0$.

We say a function $g(y):\mathbb R^{q}\rightarrow\mathbb R^{p}$ is
{\it locally bounded} if for any compact set $K\subset\mathbb
R^{q}$, there exists $M>0$ such that $\|g(y)\|\le M$ holds for all
$y\in K$.

The first item of assumption ${\mathbf A1}$ ensures that the diagonal synchronization manifold
\begin{eqnarray*}
\mathcal S=\bigg\{[x^{1\top},x^{2\top},\cdots,x^{m\top}]^{\top}\in\mathbb
R^{nm}:x^{i\top}=x^{j\top},~i,j=1,2,\cdots,m\bigg\}
\end{eqnarray*}
is an invariant manifold for (\ref{CDS}).

If $x^{1}(t)=x^{2}(t)=\cdots=x^{m}(t)=s(t)\in\mathbb R^{n}$ is the
synchronized state, then the synchronized state $s(t)$
satisfies
\begin{eqnarray}
\dot{s}(t)=f(s(t))\label{syn}.
\end{eqnarray}

Since $f(\cdot)$ is $C^{1}$-smooth, then $s(t)$ can be denoted by the corresponding
continuous semi-flow $s(t)=\vartheta^{(t)}s_{0}$ of the
intrinsic system (\ref{syn}). For $\vartheta^{(t)}$, we make following assumption:

{\bf Assumption ${\mathbf A2}$} The system (\ref{syn}) has an
asymptotically stable attractor: there exists a compact set
$A\subset R^{n}$ such that (a) $A$ is invariant through the system (\ref{syn}), i.e.,
$\vartheta^{(t)}A\subset A$ for all $t\ge 0$; (b) There exists an open bounded neighborhood $U$ of $A$ such
that $\bigcap_{t\ge 0}\vartheta^{(t)}\bar{U}=A$; (c) $A$ is topologically transitive, i.e., there exists $s_{0}\in A$
such that $\omega(s_{0})$, the $\omega$ limit set of the trajectory
$\vartheta^{(t)}s_{0}$, is equal to $A$ \cite{Ash}.
\begin{definition}
Let $A^m$ denote the Cartesian product $A\times\cdots\times A$ ($m$ times). Local complete
synchronization (synchronization for simplicity) is defined in the sense that the set $$\mathcal S\bigcap A^m=\bigg\{[x^{\top},x^{\top},\cdots,x^{\top}]^{\top}\in\mathbb
R^{nm}:x^{\top}\in A\bigg\}$$ is an asymptotically stable attractor in $\mathbb
R^{nm}$. That is, for the coupled dynamical system (\ref{CDS}), differences between components converge to
zero if the initial states are picked sufficiently near $\mathcal S\bigcap A^m$, i.e., if the components
are all close to the attractor $A$ and if their differences are sufficiently small.
\end{definition}

Next we give some lemmas which will be used later, and the proofs can be seen in Appendices.
\begin{lemma} \label{lem1} Under the assumption $\mathbf A1$, we have
\begin{eqnarray*}
\sum\limits_{j=1}^{m}\frac{\partial f^{i}}{\partial
x^{j}}(\hat{s},t)=\frac{\partial f}{\partial s}(s),
\end{eqnarray*}
for all $s\in R^{n}$ and $t\ge 0$, where $\hat{s}=[s^{\top},s^{\top},\cdots,s^{\top}]^{\top}$.
\end{lemma}
\begin{lemma}\label{lem2}Under the assumptions $\mathbf A1$ and $\mathbf A2$, there exists a
compact neighborhood $W$ of $A$ such that
$\vartheta^{(t)}W\subset\vartheta^{(t')}W$ for all $t\ge t'\ge
0$ and $\bigcap_{t\ge 0}\vartheta^{(t)}(W)=A$.
\end{lemma}

Let $\delta x(t)=[\delta
x^{1}(t)^{\top},\cdots,\delta x^{m}(t)^{\top}]^{\top}\in\mathbb
R^{nm}$, where $\delta x^{i}(t)=x^{i}(t)-s(t)\in\mathbb
R^{n}$. We have the following variational equation near the
synchronized state $s(t)$.
\begin{eqnarray*}
\delta\dot{x}^{i}(t)=\sum_{j=1}^{m}\frac{\partial
f^{i}}{\partial x^{j}}(\hat{s}(t),t)\delta x^{j}(t),i=\onetom.
\end{eqnarray*}
Or in matrix form:
\begin{eqnarray}
\delta\dot{x}(t)=DF^{t}(s(t))\delta x(t)\label{var},
\end{eqnarray}
where $DF^{t}(s(t))$ denotes the Jacobin matrix
$DF^{t}(\hat{s}(t))$ for simplicity.

From \cite{Hale}, we can give the results on the existence, uniqueness, and
continuous dependence of the equations (\ref{CDS}) and (\ref{var}).
\begin{lemma} \label{lem3} Under the assumption $\mathbf A2$, each of the
differential equations (\ref{CDS}) and (\ref{var}) has a unique solution
which is continuously dependent on the initial condition.
\end{lemma}
Thus, the solution of the linear system (\ref{var}) can be written
in matrix form.
\begin{definition}
{ Solution matrix} $U(t,t_{0},s_{0})$ of the system (\ref{var}) is
defined as follows. Let
$U(t,t_{0},s_{0})=[u^{1}(t,t_{0},s_{0}),\cdots,u^{nm}(t,t_{0},s_{0})]$,
where $u^{k}(t,t_{0},s_{0})$ denotes the $k$-th column and is the
solution of the following Cauchy problem:
\begin{eqnarray*}
\left\{\begin{array}{lll}\delta\dot{x}(t)=DF^{t}(s(t))\delta x(t)\\
s(t_{0})=s_{0}\\
\delta x(t_{0})=e^{nm}_{k}\end{array}\right..
\end{eqnarray*}
\end{definition}

Immediately, according to Lemma \ref{lem3}, we can conclude that
the solution of the following Cauchy problem
\begin{eqnarray}
\left\{\begin{array}{lll}\delta\dot{x}(t)=DF^{t}(s(t))\delta x(t)\\
s(t_{0})=s_{0}\\
\delta x(t_{0})=\delta x_{0}\end{array}\right.,\label{cauchy}
\end{eqnarray}
can be written as $\delta x(t)=U(t,t_{0},s_{0})\delta x_{0}$.

We define the time varying Jacobin matrix $DF^{t}$ by the
following way:
\begin{eqnarray*}
D\mathcal F:\mathbb R^{+}\times R^{n}&\rightarrow& 2^{\mathbb
R^{nm,nm}}\\
(t_{0},s_{0})&\mapsto& \{DF^{t}(s(t))\}_{t\ge t_{0}}
\end{eqnarray*}
with $s(t_{0})=s_{0}$, where $2^{\mathbb
R^{nm,nm}}$ is the collection of all the subsets of $\mathbb
R^{nm,nm}$.
\begin{definition} For a time varying system denoted by $D\mathcal
F$,  we can define its Hajnal diameter of the variational system (\ref{var}) as follows:
\begin{eqnarray*}
{\rm diam}(D\mathcal
F,s_{0})=\overline{\lim_{t\rightarrow\infty}}\sup_{t_{0}\ge
0}\bigg\{{{\rm diam}(U(t,t_{0},s_{0}))}\bigg\}^{\frac{1}{t}},
\end{eqnarray*}
where for a $\mathbb R^{nm,nm}$ matrix in block matrix form:
$U=(U_{ij})_{i,j=1}^{m}$ with $U_{ij}\in R^{n,n}$, its Hajnal
diameter is defined as follows:
\begin{eqnarray*}
{\rm diam}(U)=\max_{i,j}\|U_{i}-U_{j}\|,
\end{eqnarray*}
where $U_{i}=[U_{i1},U_{i2},\cdots,U_{im}]$.
\end{definition}
\begin{lemma}\label{Grounwell}(Grounwell-Beesack's inequality) If function $v(t)$ satisfies the following
condition:
\begin{eqnarray*}
v(t)\le a(t)+b(t)\int_{0}^{t}v(\tau)d\tau,
\end{eqnarray*}
where $b(t)\ge 0$ and $a(t)$ are some measurable functions, then we
have
\begin{eqnarray*}
v(t)\le
a(t)+b(t)\int_{0}^{t}a(\tau)e^{\int_{\tau}^{t}b(\theta)d\theta}d\tau,~t\ge
0.
\end{eqnarray*}
\end{lemma}

Based on the assumption $\mathbf A1$, for the solution matrix $U$, we have the following lemma:

\begin{lemma} \label{lem5}Under the assumption $\mathbf A1$, we have
\begin{enumerate}
\item
$\sum_{j=1}^{m}U_{ij}(t,t_{0},s_{0})=\breve{U}(t,t_{0},s_{0})$,
where $\breve{U}(t,t_{0},s_{0})$ denotes the solution matrix of
the following Cauchy problem:
\begin{eqnarray}
\left\{\begin{array}{ll}\dot{u}=\frac{\partial f}{\partial
s}(s(t))u\\
s(t_{0})=s_{0}\end{array}\right.;\label{syncauchy}
\end{eqnarray}
\item For any given $t\ge 0$ and the compact set $W$ given in lemma \ref{lem2},
$U(t+t_{0},t_{0},s_{0})$ is bounded for all $t_{0}\ge 0$ and
$s_{0}\in W$ and equi-continuous with respect to $s_{0}\in W$.
\end{enumerate}
\end{lemma}

Let $P=(P_{ij})_{i,j=1}^{m}$ be a
$\mathbb R^{nm,nm}$ matrix with $P_{ij}\in\mathbb R^{n,n}$
satisfying: (a) $P_{i1}=\frac{1}{\sqrt{m}}P_{0}$ for some
orthogonal matrix $P_{0}\in\mathbb R^{n,n}$ and all $i=\onetom$;
(b) $P$ is also an orthogonal matrix in $\mathbb R^{nm,nm}$. We
also write $P$ and its inverse $P^{-1}=P^{\top}$ in the form
\begin{eqnarray*}
P=[P_{1},P_{2}],\quad
P^{\top}=\left[\begin{array}{l}P^{\top}_{1}\\P_{2}^{\top}\end{array}\right],
\end{eqnarray*}
where $P_{1}=\frac{1}{\sqrt{m}}\mathbf 1_{m}\otimes P_{0}$ and
$P_{2}\in\mathbb R^{nm,n(m-1)}$. According to Lemma \ref{lem5}, we have
\begin{eqnarray*}
U(t,t_{0},s_{0})P_{1}=\frac{1}{\sqrt{m}}\mathbf
1_{m}\otimes\bigg[\breve{U}(t,t_{0},s_{0})P_{0}\bigg].
\end{eqnarray*}
Since $P_{2}^{\top}P_{1}=0$ which implies that each row of
$P_{2}^{\top}$ is located in the subspace orthogonal to the
subspace $\{\mathbf 1_{m}\otimes\xi,~\xi\in\mathbb R^{n}\}$, we
can conclude that $P_{2}^{\top}U(t,t_{0},s_{0})P_{1}=0$. Then, we
have

\begin{eqnarray*}
P^{-1}U(t,t_{0},s_{0})P=\left[\begin{array}{ll}P_{0}^{\top}\breve{U}(t,t_{0},s_{0})P_{0}&\alpha(t,t_{0},s_{0})\\
0&\tilde{U}(t,t_{0},s_{0})\end{array}\right],
\end{eqnarray*}
where $\breve{U}(t,t_{0},s_{0})$ denotes the common row sum of
$U(t,t_{0},s_{0})=(U_{ij})_{i,j=1}^{m}$ as defined in Lemma \ref{lem5},
$\tilde{U}(t,t_{0},s_{0})=P_{2}^{\top}U(t,t_{0},s_{0})P_{2}\in \mathbb R^{n(m-1),n(m-1)}$, $\alpha(t,t_{0},s_{0})\in \mathbb R^{n,n(m-1)}$ denotes a appropriate matrix we omit its accurate expression. One
can see that $\tilde{U}(t,t_{0},s_{0})$ is the solution matrix of
the following linear differential system.

\begin{definition} We define the following linear differential
system by the projection variational system of (\ref{var}) along
the directions $P_{2}$:
\begin{eqnarray}
\left\{\begin{array}{l}\dot{\phi}=D_{P}F^{t}(s(t))\phi\\
s(t_{0})=s_{0}
\end{array}\right.,\label{proj}
\end{eqnarray}
where $D_{P}F^{t}(s(t))=P_{2}^{\top}DF^{t}(s(t))P_{2}$.
\end{definition}
\begin{definition}
For any time varying variational system $D\mathcal F:\mathbb
R^{+}\times\mathbb R^{n}\rightarrow 2^{\mathbb R^{nm,nm}}$, we
define the Lyapunov exponent of the variational system (\ref{var})
as follows:
\begin{eqnarray*}
\lambda(D\mathcal
F,u,s_{0})=\overline{\lim_{t\rightarrow\infty}}\sup_{t_{0}\ge
0}\frac{1}{t}\log\|U(t,t_{0},s_{0})u\|,
\end{eqnarray*}
where $u\in\mathbb R^{nm}$ and $s(t_{0})=s_{0}$.
\end{definition}

Similarly, we can define {\it the projection Lyapunov exponents} by
the following projection time-varying variation:
\begin{eqnarray*}
D_{P}\mathcal F:\mathbb R^{+}\times\mathbb R^{n}&\rightarrow&
2^{\mathbb R^{n(m-1),n(m-1)}}\\
(t_{0},s_{0})&\mapsto&\{D_{P}F^{t}(s(t))\}_{t\ge t_{0}}
\end{eqnarray*}
i.e.,
\begin{eqnarray*}
\lambda(D_{P}\mathcal
F,\tilde{u},s_{0})=\overline{\lim_{t\rightarrow\infty}}\sup_{t_{0}\ge
0}\frac{1}{t}\log\|\tilde{U}(t,t_{0},s_{0})\tilde{u}\|,
\end{eqnarray*}
where $\tilde{u}\in\mathbb R^{n(m-1)}$ and $s(t_{0})=s_{0}$. Let
\begin{eqnarray*}
\lambda_{P}(D\mathcal F,s_{0})=\max_{\tilde{u}\in\mathbb
R^{n(m-1)}}\lambda(D_{P}\mathcal F,\tilde{u},s_{0}).
\end{eqnarray*}

Then, we have the following lemma:
\begin{lemma}\label{lem6}
$\lambda_{P}(D\mathcal F,s_{0})=\log{\rm diam}(D\mathcal
F,s_{0})$.
\end{lemma}

\begin{remark}
From Lemma \ref{lem6}, we can see that the
largest projection Lyapunov exponent is independent of the choice
of matrix $P$.
\end{remark}

Consider the time-varying driven by some metric dynamical
system (MDS) $(\Omega,\mathcal B,\mathbb P,\varrho^{(t)})$, where $\Omega$
is the compact state space, $\mathcal B$ is the $\sigma$-algebra,
$\mathbb P$ is the probability measure, and $\varrho^{(t)}$ is a
continuous semi-flow. Then, the variational equation (\ref{var}) is
independent of the initial time $t_{0}$ and can be rewritten as
follows:
\begin{eqnarray}
\left\{\begin{array}{l}\dot{\phi}=DF(s(t),\varrho^{(t)}\omega_{0})\phi\\
s(0)=s_{0}\end{array}\right.\label{varrds}
\end{eqnarray}
In this case, we denote {\it the solution matrix, the projection
solution matrix}, and {\it the solution matrix} on the synchronization
space by $U(t,s_{0},\omega_{0})$, $\tilde{U}(t,s_{0},\omega_{0})$,
and $\breve{U}(t,s_{0},\omega_{0})$, respectively. For simplicity,
we write them as $U(t)$, $\tilde{U}(t)$, and $\breve{U}(t)$,
respectively. Also, we write {\it the Lyapunov exponents} and {\it the projection Lyapunov exponent} as follows:
\begin{align*}
\lambda(D\mathcal F, u, s_{0},\omega_{0})&=\overline{\lim\limits_{t\rightarrow\infty}}\frac{1}{t}\log\|U(t,s_{0},\omega_{0})u\|,~~
\lambda(D\mathcal F,s_{0},
\omega_{0})=\max\limits_{u\in\mathbb R^{nm}}\lambda(D\mathcal F,u,s_{0},\omega_{0}),\\
\lambda_{P}(D\mathcal F,u, s_{0},\omega_{0})&=\overline{\lim\limits_{t\rightarrow\infty}}\frac{1}{t}\log\|\tilde{U}(t,s_{0},\omega_{0})u\|,~~
\lambda_{P}(D\mathcal F,s_{0},
\omega_{0})=\max\limits_{u\in\mathbb R^{n(m-1)}}\lambda(D_{P}\mathcal F,\tilde{u},s_{0},\omega_{0}).\\
\end{align*}
We add the following assumption:

{\bf Assumption ${\mathbf
A3}$} (a) $\varrho^{(t)}$ is a continuous semi-flow; (b) $DF(s,\omega)$ is a continuous map for all $(s,\omega)\in
\mathbb R^{n}\times \Omega$.

The following are involving linea differential systems. For more
details, we refer the readers to \cite{Adr}.  For a continuous
scalar function $u(t)$, we denote its {\it Lyapunov exponent} by
\begin{eqnarray*}
\chi[u(t)]=\overline{\lim\limits_{t\rightarrow\infty}}\frac{1}{t}\log|u(t)|.
\end{eqnarray*}
The following properties will be used later:
\begin{enumerate}
\item
$\chi[\prod\limits_{k=1}^{n}c_{k}u_{k}(t)]\le\sum\limits_{k=1}^{n}\chi[u_{k}(t)]$,
where $c_{k}$, $k=\oneton$, are constants;

\item if
$\lim\limits_{t\rightarrow\infty}\frac{1}{t}\log|u(t)|=\alpha$,
which is finite, then $\chi[\frac{1}{u(t)}]=-\alpha$;

\item $\chi[u(t)+v(t)]\le\max\{\chi[u(t)],\chi[v(t)]\}$;

\item for a vector-value or matrix-value function $U(t)$, we
define $\chi[U(t)]=\chi[\|U(t)\|]$.
\end{enumerate}

For the following linear differential system:
\begin{eqnarray}
\dot{x}(t)=A(t)x(t)\label{a}
\end{eqnarray}
where $x(t)\in\mathbb R^{n}$, a transformation $x(t)=L(t)y(t)$ is
said to be a {\it Lyapunov transformation} if $L(t)$ satisfies :
\begin{enumerate}
\item $L(t)\in C^{1}[0,+\infty)$;

\item $L(t)$, $\dot{L}(t)$, $L^{-1}(t)$ are bounded for all $t\ge
0$.
\end{enumerate}
It can be seen that the class of Lyapunov transformations form a
group and the linear system for $y(t)$ should be
\begin{eqnarray}
\dot{y}(t)=B(t)y(t),\label{b}
\end{eqnarray}
where $B(t)=L^{-1}(t)A(t)L(t)-L^{-1}(t)\dot{L}(t)$. Then, we say
system (\ref{b}) is a {\it reducible system } of system (\ref{a}).
We define the {\it adjoint system} of (\ref{a}) by
\begin{eqnarray}
\dot{x}(t)=-A^{*}(t)x(t).\label{c}
\end{eqnarray}
If letting $V(t)$ be the fundamental matrix of (\ref{a}), then
$[V^{-1}(t)]^{*}$ is the fundamental matrix of (\ref{c}). Thus, we
say the system (\ref{a}) is a {\it regular system} if the adjoint
systems (\ref{a}) and (\ref{c}) have convergent Lyapunov exponent series:
$\{\alpha_{1},\cdots,\alpha_{n}\}$ and
$\{\beta_{1},\cdots,\beta_{n}\}$, respectively, which satisfy
$\alpha_{i}+\beta_{i}=0$ for $i=\oneton$; or its reducible system
(\ref{b}) is also regular.

\begin{lemma}\label{lem7}Suppose that the assumptions $\mathbf A1$, $\mathbf A2$, and $\mathbf A3$ are
satisfied. Let
$\{\sigma_{1},\sigma_{2},\cdots,\sigma_{n},\sigma_{n+1},\cdots,\sigma_{nm}\}$
be the Lyapunov exponents of the variational system
(\ref{varrds}), where $\{\sigma_{1},\cdots,\sigma_{n}\}$
correspond to the synchronization space and the remaining
correspond to the transverse space. Let $\lambda_{T}(D\mathcal
F,s_{0},\omega_{0})=\max\limits_{i\ge n+1}\sigma_{i}$ and
$\lambda_{S}(D\mathcal F,s_{0},\omega_{0})=\max\limits_{1\le i\le
n}\sigma_{i}$. If
(a) the linear system (\ref{syncauchy}) is a regular system,
(b) $\|DF(s(t),\varrho^{(t)}\omega_{0})\|\le M$ for all $t\ge0$,
(c) $\lambda_{P}(D\mathcal
F,s_{0},\omega_{0})\ne\lambda_{S}(D\mathcal F,s_{0},\omega_{0})$,
then $\lambda_{T}(D\mathcal
F,s_{0},\omega_{0})=\lambda_{P}(D\mathcal F,s_{0},\omega_{0})$.
\end{lemma}

\section{General synchronization analysis}
In this section we provide a methodology based on the previous theoretical analysis to judge whether a general differential system can be synchronized or not.
\begin{theorem}\label{thm1}
Suppose that $W\in \mathbb R^{n}$ is the compact subset given in Lemma \ref{lem2},
and assumptions $\mathbf A1$ and $\mathbf A2$ are satisfied. If
\begin{eqnarray}
\sup_{s_{0}\in W} {\rm diam}(D\mathcal F,s_{0})<1,\label{cond1}
\end{eqnarray}
then the coupled system (\ref{CDS}) is synchronized.
\end{theorem}
{\bf Proof}: The main techniques of the proof come from \cite{Ash,LAJ} with some modifications.
Let $\vartheta^{(t)}$ be the semi-flow of the
uncoupled system (\ref{syn}). By the condition (\ref{cond1}), there exist $d$
satisfying $\sup_{s_{0}\in W} {\rm diam}(D\mathcal
F,s_{0})<d<1$ and $T_{1}\ge 0$ such that
$d^{T_{1}}<\frac{1}{3}$, and $r_{0}=\inf\bigg\{r>0,~\mathcal
O(\vartheta^{(T_{1})} W,r)\subset W\bigg\}>0$. For each $s_{0}\in
W$, there must exists $t(s_{0})\ge T_{1}$ such that ${\rm
diam}(U(t_{0}+t(s_{0}),t_{0},s_{0}))<d^{t(s_{0})}$ for all
$t_{0}\ge 0$. According to the equi-continuity of
$U(t_{0}+t(s_{0}),t_{0},s_{0})$, there exists $\delta>0$ such that
for any $s_{0}^{'}\in\mathcal O(s_{0},\delta)$, ${\rm
diam}(U(t_{0}+t(s_{0}),t_{0},s_{0}^{'}))<d^{t(s_{0})}$ for all
$t_{0}\ge 0$. According to the compactness of $ W$, there exists a
finite positive number set $\mathcal
T=\{t_{1},t_{2},\cdots,t_{v}\}$ with $t_{j}\ge T_{1}$ for all
$j=1,2,\cdots,v$ such that for any $s_{0}\in W$, there exists
$t_{j}\in\mathcal T$ such that ${\rm
diam}(U(t_{0}+t_{j},t_{0},s_{0}))<\frac{1}{3}$ for all
$t_{0}\ge 0$. Let $x(t)$ be the collective  states
$\{x^{1}(t),\cdots,x^{m}(t)\}$ which is the solution of the
coupled system (\ref{CDS}) with initial condition
$x^{i}(t_{0})=x^{i}_{0}$, $i=\onetom$. And, let $s(t)$ be the solution
of the synchronization state equation (\ref{syn}) with initial
condition
$s(t_{0})=\bar{x}_{0}=\frac{1}{m}\sum_{j=1}^{m}x^{j}_{0}\in W$.
Then, letting $\Delta x^{i}(t)=x^{i}(t)-s(t)$, we have
\begin{align*}
\Delta
\dot{x}_{k}^{i}(t)=f_{k}^{i}(x^{1}(t),\cdots,x^{m}(t),t)-f_{k}(s(t))=\sum_{j=1}^{m}\sum_{l=1}^{n}\frac{\partial
f^{i}_{k}}{\partial x^{j}_{l}}(\xi^{ij}_{kl}(t),t)\Delta
x^{j}_{l}(t),
\end{align*}
where $\xi^{ij}_{kl}(t)\in \mathbb R^{mn}$, $i,j=\onetom$,
$k,l=\oneton$, are obtained by the mean value principle of the differential functions. Letting
$DF^{t}(\xi(t))=(\frac{\partial f^{i}_{k}}{\partial
x^{j}_{l}}(\xi^{ij}_{kl}(t),t))$, we can write the equations above
in matrix form:
\begin{eqnarray*}
\Delta \dot{x}(t)=DF^{t}(\xi(t))\Delta x(t),
\end{eqnarray*}
and denote its solution matrix by
$\hat{U}(t+t_{0},t_{0},x_{0})=(\hat{U}_{ij}(t+t_{0},t_{0},x_{0}))_{i,j=1}^{m}$.
Then, for any $t>0$ there exists $K_{2}>0$ such that
$\|DF^{t+t_{0}}(\xi(t+t_{0}))\|\le K_{2}$ for all
$t\in\mathcal T$ and $t_{0}\ge 0$ according to the $3$-th item of the
assumption $\mathbf A1$. Then, we have
\begin{align*}
\Delta
x^{i}_{k}(t+t_{0})
=x^{i}_{0}-\bar{x}_{0k}
+\int_{t_{0}}^{t+t_{0}}\sum_{j=1}^{m}\sum_{l=1}^{n}\frac{\partial
f^{i}_{k}}{\partial x^{j}_{l}}(\xi^{ij}_{kl}(\tau),\tau)\Delta
x^{j}_{l}(\tau)d\tau,\\
\sum_{j=1}^{m}\sum_{k=1}^{n}\|\Delta
x^{i}_{k}(t+t_{0})\|\le\sum_{j=1}^{m}\sum_{k=1}^{n}\|x^{i}_{0k}-\bar{x}_{0k}\|
+K_{2}\int_{t_{0}}^{t+t_{0}}\sum_{j=1}^{m}\sum_{l=1}^{n}\|\Delta
x^{j}_{l}(\tau)\|d\tau.
\end{align*}
By Lemma \ref{Grounwell}, we have
\begin{eqnarray*}
\sum_{j=1}^{m}\sum_{l=1}^{n}\|\Delta
x^{j}_{l}(t+t_{0}) \|\le
e^{K_{2}t}\sum_{j=1}^{m}\sum_{l=1}^{n}\|x^{j}_{0l}-\bar{x}_{0l}\|.
\end{eqnarray*}
Let
\begin{eqnarray*}
W_{\alpha}=\{x=[{x^{1}}^{\top},\cdots,{x^{m}}^{\top}]^{\top}:~\bar{x}\in
W,\sum_{j=1}^{m}\|x^{j}-\bar{x}\|\le\alpha\}.
\end{eqnarray*}
Picking $\alpha$ sufficiently small such that for each $x_{0}\in
W_{\alpha}$, there exists $t\in\mathcal T$ such that
$\sum_{j=1}^{m}\|\Delta
x^{j}(t+t_{0})\|<\frac{r_{0}}{2}$ and ${\rm
diam}(\hat{U}(t+t_{0},t_{0},x_{0}))<\frac{1}{2}$for all
$t_{0}\ge 0$.

Thus, we are to prove synchronization step by step.

For any $x_{0}\in
W_{\alpha}$, there exists $t'=t(x_{0})\in\mathcal T$ such that
\begin{align*}
\|x^{i}(t'+t_{0})-x^{j}(t'+t_{0})\|&=\|\Delta
x^{i}(t'+t_{0})-\Delta
x^{j}(t'+t_{0})\|\\
&\le\sum_{k=1}^{m}\|\hat{U}_{ik}(t'+t_{0},t_{0},x_{0})-\hat{U}_{jk}(t'+t_{0},t_{0},x_{0})\|\|\Delta
x^{k}_{0}\|\\
&\le{\rm
diam}(\hat{U}(t'+t_{0},t_{0},x_{0}))\max_{i,j}\|x^{i}_{0}-x^{j}_{0}\|\\
&\le\frac{1}{2}\max_{i,j}\|x^{i}_{0}-x^{j}_{0}\|
\end{align*}
Therefore, we have
$\max_{i,j}\|x^{i}(t'+t_{0})-x^{j}(t'+t_{0})\|\le\frac{1}{2}\max_{i,j}\|x^{i}_{0}-x^{j}_{0}\|$,
which implies that $\bar{x}(t'+t_{0})\in W$ and $x(t'+t_{0})\in
W_{\frac{\alpha}{2}}$.

Then, re-initiated with time $t'+t_{0}$ and condition
$x(t'+t_{0})$, continuing with the phase above, we can obtain that
$\lim_{t\rightarrow\infty}\max_{i,j}\|x^{i}(t)-x^{j}(t)\|=0$.
Namely, the coupled system (\ref{CDS}) is synchronized.
Furthermore, from the proof, we can conclude that the convergence
is exponential with rate $O(\delta^{t})$ where
$\delta=\sup_{s_{0}\in W}{\rm diam}(D\mathcal
F^{t},s_{0})$, and uniform with respect to $t_{0}\ge 0$ and
$x_{0}\in W_{\alpha}$. This completes the proof.

\begin{remark}
According to the assumption $\mathbf A2$ that attractor $A$ is asymptotically stable and the
properties of the compact neighbor $W$ given in Lemma \ref{lem2}, we can conclude
that the quantity $$\sup_{s_{0}\in W}{\rm diam}(D\mathcal
F,s_{0})$$ is independent on the choice of $W$.
\end{remark}

If the time-variation is driven by some MDS $\{\Omega,\mathcal
B,\mathbb P P,\varrho^{(t)}\}$ and there exists a metric dynamical system
$\{ W\times\Omega,{\bf F},{\bf P},\pi^{(t)}\}$, where $\bf F$ is
the product $\sigma-$ algebra on $ W\times\Omega$, $\bf P$ is the
probability measure, and
$\pi^{(t)}(s_{0},\omega)=(\theta^{(t)}s_{0},\varrho^{(t)}\omega)$.
From Theorem 1, we have
\begin{corollary}
Suppose that the conditions in Lemma\ref{lem7} are satisfied, $ W\times \Omega$ is compact in the topology
defined in this MDS, the semi-flow $\pi^{(t)}$ is continuous, and on $ W\times\Omega$
the Jacobian matrix $DF(\theta^{(t)}s_{0},\varrho^{(t)}\omega)$ is
continuous. Let $\{\sigma_{i}\}_{i=1}^{nm}$ be
the Lyapunov exponents of this MDS with multiplicity and
$\{\sigma_{i}\}_{i=1}^{n}$ correspond the synchronization space.
If
\begin{eqnarray*}
\sup\limits_{{\bf P}\in Erg_{\pi}(
W\times\Omega)}\sup\limits_{i\ge n+1}\sigma_{i}<0,
\end{eqnarray*}
where $Erg_{\pi}( W\times\Omega)$ denotes the ergodic probability
measure set supported in the MDS $\{ W\times\Omega,{\bf F},{\bf
P},\pi^{(t)}\}$, then the coupled system (\ref{CDS}) is
synchronized.
\end{corollary}

\section{Synchronization of LCODEs with identity inner coupling matrix and time-varying couplings}

In this section we study synchronization in
linearly coupled ordinary differential equation systems (LCODEs) with time-varying couplings. Considering the following LCODEs  with identity inner coupling matrix:
\begin{eqnarray}
\dot{x}^{i}(t)=f(x^{i}(t))+\sigma\sum_{j=1}^{m}l_{ij}(t)x^{j}(t),~i=\onetom,\label{LCODE}
\end{eqnarray}
where $x^{i}(t)\in\mathbb R^{n}$ denotes the state variable of the
$i$-th node, $f(\cdot):\mathbb R^{n}\rightarrow\mathbb R^{n}$ is a
differential map, $\sigma\in \mathbb R^{+}$ denotes coupling strength, and $l_{ij}(t)$ denotes the coupling coefficient from node $j$ to $i$ at time $t$, for all $i\ne j$, are supposed to satisfy the
following assumption. Here, we highlight that the inner coupling matrix is the identity matrix.

{\bf Assumption $\mathbf A4$} (a)
$l_{ij}(t)\ge 0$, $i\ne j$, are measurable and
$l_{ii}(t)=-\sum_{j=1,j\ne i}^{m}l_{ij}(t)$; (b) There exists $M_{1}>0$ such that $|l_{ij}(t)|\le M_{1}$ for all $i,j=\onetom$.

In the following we omit coupling strength $\sigma$ since it is just a positive scalar. Similarly, we can define the Hajnal diameter of the following
linear system:
\begin{eqnarray}
\dot{u}(t)= \sigma L(t)u(t).\label{linear}
\end{eqnarray}
Let $V(t)=(v_{ij}(t))_{i,j=1}^{m}$ be the fundamental solution
matrix of the system (\ref{linear}). Then, its solution matrix can
be written as $V(t,t_{0})=V(t){V}(t_{0})^{-1}$. Thus, the Hajnal
diameter of the system (\ref{linear}) can be defined as follows:
\begin{eqnarray*}
{\rm diam}( \mathcal
L)=\overline{\lim_{t\rightarrow\infty}}\sup_{t_{0}\ge
0}\bigg[{\rm diam}(V(t,t_{0}))\bigg]^{1/t}.
\end{eqnarray*}

 By Theorem \ref{thm1}, we have
\begin{theorem}
Suppose the assumptions $\mathbf A1$, $\mathbf A2$, and $\mathbf A4$ are satisfied. Let
$\mu$ be the largest Lyapunov exponent of the synchronized system
$\dot{s}(t)=f(s(t))$, i.e.,
\begin{eqnarray*}
\mu=\sup_{s_{0}\in W}\max_{u\in\mathbb
R^{n}}\lambda(D f,u,s_{0}).
\end{eqnarray*}
If $\log({\rm diam}(\mathcal L))+\mu<0$, then the LCODEs
(\ref{LCODE}) is synchronized.
\end{theorem}
{\bf Proof}: Considering the variational equation of (\ref{LCODE}):
\begin{eqnarray}
\delta\dot{x}(t)=\bigg[I_{m}\otimes Df(s(t))+\sigma L(t)\otimes
I_{n}\bigg]\delta x(t).\label{varLCODE}
\end{eqnarray}
Let $\breve{U}(t,t_{0},s_{0})$ be the solution matrix of the
synchronized state system (\ref{syncauchy}) and
$V(t,t_{0})=(v_{ij}(t,t_{0}))_{i,j=1}^{m}$ be the solution matrix of the linear
system (\ref{linear}). We can see that
$V(t,t_{0})\otimes\breve{U}(t,t_{0},s_{0})$ is the solution matrix
of the variational system (\ref{varLCODE}). Then,
\begin{align*}
&{\rm
diam}(V(t+t_{0},t_{0})\otimes\breve{U}(t+t_{0},t_{0},s_{0}))\\
&=\max_{i,j=1,\cdots,m}\sum_{k=1}^{m}
|v_{ik}(t+t_{0},t_{0})-v_{jk}(t+t_{0},t_{0})|\times\|\breve{U}(t+t_{0},t_{0},s_{0})\|\\
&={\rm diam}(V(t+t_{0},t_{0}))\|\breve{U}(t+t_{0},t_{0},s_{0})\|.
\end{align*}
This implies that the Hajnal diameter of the variational system
(\ref{varLCODE}) is less than $e^{\mu}{\rm diam}(\mathcal L)$.
This completes the proof according to Theorem 1.

For the linear system (\ref{linear}), we firstly have
\begin{lemma}\cite{Liubo}\label{lem8}
$V(t,t_{0})$ is a stochastic matrix.
\end{lemma}

From Lemmas \ref{lem6} and
\ref{lem7}, we have
\begin{corollary}
$\log{\rm diam}(\mathcal L)=\lambda_{P}(\mathcal L)$, where $\lambda_{P}(\mathcal L)$ denotes the largest one of all the projection Lyapunov exponents of system (\ref{LCODE}). Moreover, if
the conditions in Lemma 7 are satisfied, then $\log{\rm
diam}(\mathcal L)=\lambda_{T}(\mathcal L)$, where $\lambda_{T}(\mathcal L)$ denotes the largest one of all the Lyapunov exponents corresponding to the transverse space, i.e., the space orthogonal to the synchronization space.
\end{corollary}

If $L(t)$ is periodic, we have
\begin{corollary}
Suppose that $L(t)$ is periodic. Let $\varsigma_{i}$, $i=\onetom$,
are the Floquet multipliers of the linear system (\ref{linear}).
Then, there exists one multiplier denoted by $\varsigma_{1}=1$ and
${\rm diam}(\mathcal L)=\max\limits_{i\ge 2}\varsigma_{i}$.
\end{corollary}

If $L(t)=L(\varrho^{(t)}\omega)$ is driven by some MDS
$\{\Omega,\mathcal B,P,\varrho^{(t)}\}$, from Corollaries 1 and 2,
we have
\begin{corollary}
Suppose $L(\omega)$ is continuous on $\Omega$ and conditions in
Lemma 7 are satisfied. Let $\mu=\sup\limits_{s_{0}\in
W}\max\limits_{u\in \mathbb R^{n}}\lambda(Df,u,s_{0})$,
$\varsigma_{i}$, $i=\onetom$, be the Lyapunov exponents of the
linear system (\ref{linear}) with $\varsigma_{1}=0$, and
$\varsigma=\sup\limits_{P\in
Erf_{\theta}(\Omega)}\max\limits_{i\ge 2}\varsigma_{i}$. If
$\mu+\varsigma<0$, then the coupled system (\ref{LCODE}) is
synchronized.
\end{corollary}



Let $\mathcal I$ be the set consisting of
all compact time intervals in $[0,+\infty)$ and $\mathcal G$ be
the the set consisting of all graph with vertex set
$N=\{1,2,\cdots,m\}$.

Define
\begin{eqnarray*}
G:\mathcal I\times R^{+}&\rightarrow& \mathcal G\\
(I=[t_{1},t_{2}],\delta)&\mapsto& G(I,\delta)
\end{eqnarray*}
where $G(I,\delta)=\{N,E\}$ is a graph with vertex set $N$ and its
edge set $E$ is defined as follows: there exists an edge from
vertex $j$ to vertex $i$  if and only if
$\int_{t_{1}}^{t_{2}}l_{ij}(\tau)d\tau>\delta$. Namely, we say
that {\it there is a $\delta$-edge from vertex $j$ to $i$ across
$I=[t_{1},t_{2}]$}.
\begin{definition}
We say that the LCODEs (\ref{LCODE}) has a
$\delta$-spanning tree across the time interval $I$ if the corresponding graph
$G(I,\delta)$ has a spanning tree.
\end{definition}
For a stochastic matrix
$V=(v_{ij})_{i,j=1}^{m}$, let
\begin{eqnarray*}
\eta(V)=\min_{i,j}\|v_{i}\wedge v_{j}\|_{1}
\end{eqnarray*}
where $v_{i}=[v_{i1},\cdots,v_{im}]$, $i=\onetom$, and
$v_{i}\wedge
v_{j}=[\min\{v_{i1},v_{j1}\},\cdots,\min\{v_{im},v_{jm}\}]^{\top}$.
Then, we can also define that
$V$ is {\it $\delta$-scrambling} if $\eta(V)>\delta$.

\begin{theorem}
Suppose the assumption $\mathbf A4$ is satisfied. ${\rm diam}(\mathcal L)<1$ if and only if there exist $\delta>0$
and $T>0$ such that the LCODEs (\ref{LCODE}) has a $\delta$-spanning tree across any $T$-length time
interval.
\end{theorem}
\begin{remark}
Different from \cite{Sti}, we don't need to assume that $L(t)$ has zero column sums and  the time-average is achieved sufficiently fast.
\end{remark}
Before proving this theorem, we need the following lemma:
\begin{lemma}\label{lem9}
If the LCODEs (\ref{LCODE}) has a $\delta$-
spanning tree across any $T$-length time interval, then there
exist $\delta_{1}>0$ and $T_{1}>0$ such that
$V(t,t_{0})$ is $\delta_{1}$-scrambling for any
$T_{1}$-length time interval.
\end{lemma}

{\bf Proof of Theorem 3}: Sufficiency. From Lemma \ref{lem9}, we
can conclude that there exist $\delta_{1}>0$, $\delta'>0$, and
$T_{1}>0$ such that $V(t,t_{0})$ is $\delta_{1}$-scrambling across
any $T_{1}$-length time interval and  $\inf_{t_{0}\ge
0}\eta(V(T_{1}+t_{0},t_{0}))>\delta'$. For any $t\ge t_{0}$, let $t-t_{0}=pT_{1}+T'$, where $p$ is an
integer and $0\le T'< T_{1}$ and $t_{l}=t_{0}+lT_{1}$, $0\le l\le
p$. Then, we have
\begin{align*}
{\rm diam}(V(t,t_{0}))&={\rm
diam}(V(t,t_{p})\prod_{l=1}^{p}V(t_{l},t_{l-1}))\le {\rm diam}(\prod_{l=1}^{p}V(t_{l},t_{l-1}))\\
&\le2\prod_{l=1}^{p}(1-\eta(V(t_{l},t_{l-1})))\le
2(1-\delta')^{\lfloor \frac{t-t_{0}}{T_{1}}\rfloor}.
\end{align*}
(For the first inequality, we use the results in \cite{Haj1,Haj2}.) This implies ${\rm diam}(\mathcal L)\le
(1-\delta')^{\frac{1}{T_{1}}}<1$.

Necessity. Suppose that for any $T\ge 0$ and
$\delta>0$, there exists $t_{0}=t_{0}(T,\delta)$,
$\int_{t_{0}}^{T+t_{0}}L(\tau)d\tau$ does not have a
$\delta$-spanning tree. According to the condition, there exist
$1>d>{\rm diam}(\mathcal L)$, $\epsilon>0$, and $T'>0$ such that
${\rm diam}(V(t+t_{0}))<d^{t}$ for all $t_{0}\ge 0$ and
$t\ge T'$ and $d^{T'}<1-\epsilon$.
Thus, picking $T>T'$, $\delta=m^{-3}e^{-M_{1}mT}\epsilon/2$,
$t_{1}=t_{0}(T,\delta)$, and
$L^{'}=(l^{'}_{ij})_{i,j=1}^{m}=(\int_{T}^{T+t_{1}}l_{ij}(\tau)d\tau)_{i,j=1}^{m}$,
there exist two vertex set $J_{1}$ and $J_{2}$ such that $
l^{'}_{ij}\le \delta$ if $i\in J_{1}$ and $j\notin J_{1}$, or
$i\in J_{2}$ and $j\notin J_{2}$.  For each $i\in J_{1}$ and
$j\notin J_{1}$, we have
\begin{align*}
\dot{v}_{ij}(t)=l_{ii}(t)v_{ij}(t)+\sum_{k\in J_{1}}^{k\ne
i}l_{ik}(t)v_{kj}(t)+\sum_{k\notin
J_{1}}l_{ik}(t)v_{kj}(t)\le M_{1}\sum_{k\in J_{1}}^{k\ne
i}v_{kj}(t)+\sum_{k\notin J_{1}}l_{ik}(t).
\end{align*}
Then,
\begin{align*}
\sum_{i\in J_{1},j\notin
J_{1}}\dot{v}_{ij}(t)&\le M_{1}\sum_{i\in J_{1},k\in
J_{1}}^{k\ne i,j\notin J_{1}}v_{kj}(t)+\sum_{i\in J_{1},k\notin
J_{1}}^{j\notin J_{1}}l_{ik}(t)\\
&=M_{1}(\# J_{1}-1)\sum_{k\in J_{1}}^{j\notin
J_{1}}v_{kj}(t)+(m-\# J_{1})\sum_{i\in J_{1}}^{k\notin
J_{1}}l_{ik}(t).
\end{align*}
Let $v(t)=\sum_{i\in J_{1},j\notin J_{1}}v_{ij}(t)$.
According to Lemma \ref{Grounwell}, we have
\begin{align*}
v(T+t_{1})&\le e^{M_{1}(\# J_{1}-1)T}(m-\#
J_{1})\int_{t_{1}}^{T+t_{1}}\sum_{i\in J_{1}}^{j\notin
J_{1}}l_{ij}(\tau)d\tau\\
&\le (m-\# J_{1})e^{M_{1}(\#J_{1}-1)T}\# J_{1}(m-\#
J_{1})\delta\le m^{3}e^{mM_{1}T}\delta\le \frac{\epsilon}{2}.
\end{align*}
Similarly, we can conclude that $\sum_{i\in J_{l},j\notin
J_{l}}v_{ij}(T+t_{1})\le \frac{\epsilon}{2}$ for all
$l=1,2$. Without loss of generality, we suppose
$J_{1}=\{1,2,\cdots,p\}$ and $J_{2}=\{p+1,p+2,\cdots,p+q\}$, where
$p$ and $q$ are integers with $p+q\le m$. Then, we can write
$V(T+t_{1},t_{1})$ in the following matrix form:
\begin{eqnarray*}
V(T+t_{1},t_{1})=\left[\begin{array}{lll}X_{11}&X_{12}&X_{13}\\
X_{21}&X_{22}&X_{23}\\
X_{31}&X_{32}&X_{33}\end{array}\right],
\end{eqnarray*}
where $X_{11}\in\mathbb R^{p,p}$ and $X_{22}\in R^{q,q}$
correspond the vertex subset $J_{1}$ and $J_{2}$ respectively.
Immediately, we have
$\|X_{12}\|_{\infty}+\|X_{13}\|_{\infty}+\|X_{21}\|_{\infty}+\|X_{23}\|_{\infty}\le\epsilon$. Let $v=\left[\begin{array}{l}{\bf
1}_{p}\\0\\0\end{array}\right]$. We let
\begin{eqnarray*}
V(t_{1}+T,t_{1})v=\left[\begin{array}{l}X_{11}{\bf
1}_{p}\\
X_{21}{\bf
1}_{p}\\
X_{31}{\bf 1}_{p}\end{array}\right].
\end{eqnarray*}
Let
$u=\left[\begin{array}{l}u^{1}\\u^{2}\\u^{3}\end{array}\right]=[u_{1},\cdots,u_{m}]^{\top}$
with $u^{i}=[u^{i}_{1},\cdots,u^{i}_{p_{i}}]^{\top}=X_{i1}{\bf
1}_{p}$ and  $p_{1}=p$, $p_{2}=q$, $p_{3}=m-p-q$. Then,
\begin{eqnarray*}
\max_{i,j}|u_{i}-u_{j}|\ge\max_{k,l}|u^{1}_{k}-u^{2}_{j}|
\ge1-\|X_{12}\|_{\infty}-\|X_{13}\|_{\infty}-\|X_{21}\|_{\infty}-\|X_{23}\|_{\infty}\ge 1-\epsilon.
\end{eqnarray*}
Also,
\begin{eqnarray*}
\max_{i,j}|u_{i}-u_{j}|\le{\rm diam}(V(t_{1}+T,t_{1}))\le
d^{T}.
\end{eqnarray*}
This implies $d^{T}\ge 1-\epsilon$ which leads contradiction with
$d^{T}<1-\epsilon$. Therefore, we can conclude the necessity.

\section{Consensus analysis of multi-agent system with directed time-varying graphs}
If let $n=1$, $f\equiv 0$ and $\sigma=1$ in system (\ref{LCODE}), then we have
\begin{eqnarray}\label{consen}
\dot{x}^{i}(t)=\sum_{j=1}^{m}l_{ij}(t)x^{j}(t),~i=\onetom.
\end{eqnarray}
In this case, if the assumption $\mathbf A4$ is satisfied, then the synchronization analysis of system (\ref{consen}) become another important research field named consensus problems.
\begin{definition}
We say the differential system (\ref{consen}) reaches consensus if for any $x(t_0)\in\mathbb R^{m}$, $\|x^{i}(t)-x^{j}(t)\|\rightarrow 0$ as $t\rightarrow \infty$ for all $i,j\in \{1,\cdots,m\}$.
\end{definition}

In graph view, the coefficients matrix of (\ref{consen}) $L(t)=(l_{ij}(t))\in \mathbb R^{m,m}$ is equal to the negative {\it graph Laplacian} associated with the digraph $G(t)$ at time $t$, where $G(t)=(\mathcal V, \mathcal E(t), \mathcal A(t))$ is a weighted digraph (or directed graph) with $m$ vertices, the set of nodes $\mathcal V =\{v_1,\cdots,v_m\}$, set of edges $\mathcal E(t) \subseteq \mathcal V \times \mathcal V$, and the weighted adjacency matrix $\mathcal A(t) =(a_{ij}(t))$ with nonnegative adjacency elements $a_{ij}(t)$. An edge of $G(t)$ is denoted by $e_{ij}(t)=(v_i, v_j)\in \mathcal E(t)$ if there is an directed edge from vertex $i$ to vertex $j$ at time $t$. The adjacency elements associated with the edges of the graph are positive, i.e., $e_{ij}(t)\in \mathcal E(t)\iff a_{ij}(t)>0$, for all $i, j\in\{1, 2,\cdots, m\}$. It is assumed that $a_{ii}(t)
=0$ for all $i\in \{1, 2,\cdots, m\}$.
The in-degree and out-degree of node $v_i$ at time $t$ are, respectively, defined as follows:
\begin{eqnarray*}
deg_{in}(v_i(t))=\sum\limits_{j=1}^{N}a_{ji(t)},~~~deg_{out}(v_i(t))=\sum\limits_{j=1}^{N}a_{ij(t)}.
\end{eqnarray*}
The degree matrix of digraph $G(t)$ at time $t$ is defined as $D(t)=diag(deg_{out}(v_1(t)), \cdots, deg_{out}(v_m(t)))$. The {\it graph Laplacian} associated with the digraph $G(t)$ at time $t$ is defined as
\begin{eqnarray}
-L(t)=\mathcal L (G(t))=D(t)-\mathcal A(t).
\end{eqnarray}

Let $G(I,\delta)$ defined as before. We say that the digraph $G(t)$ has a
$\delta$-spanning tree across the time interval $I$ if
$G(I,\delta)$ has a spanning.

\begin{theorem}
Suppose the assumption $\mathbf A4$ is satisfied. The system (\ref{consen}) reaches consensus if and only if there exist $\delta>0$
and $T>0$ such that the corresponding digraph $G(t)$ has a $\delta$-spanning tree across any $T$-length time
interval.
\end{theorem}
{\bf Proof}: Since $f\equiv 0$, we have $\mu=0$ in Theorem 2. This completes the proof according to Theorem 3 and Theorem 2.

\begin{remark}
This theorem is a part of the Theorem 1 in \cite{Mor}.
\end{remark}

\section{Numerical examples}

\begin{figure}[hbt]
\centering
\includegraphics[width=.6\textwidth]{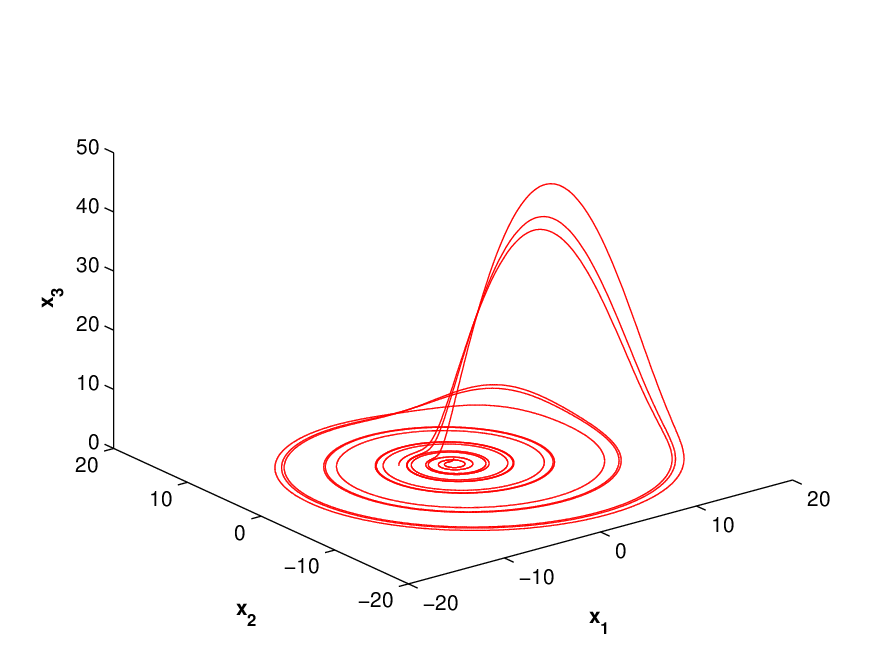}
\caption{The dynamical behavior of the R\"{o}ssler system (\ref{ross}) with $a=0.165, b=0.2$, and $c=10$.}
\label{fig:1}
\end{figure}

\begin{figure}[hbt]
\centering
\includegraphics[width=2.5in]{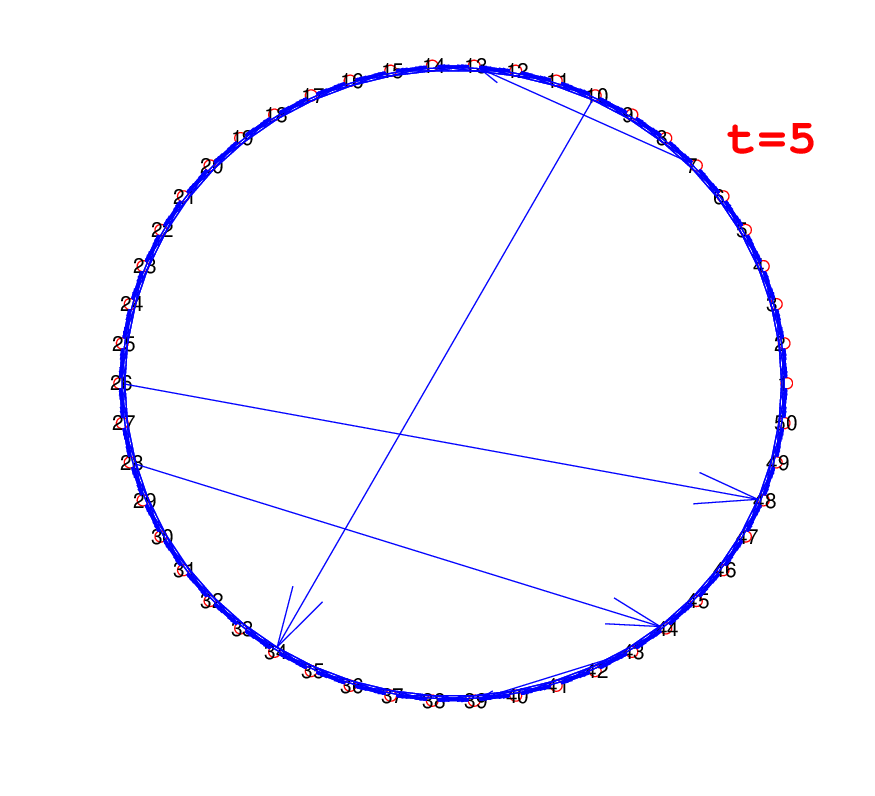}
\includegraphics[width=2.5in]{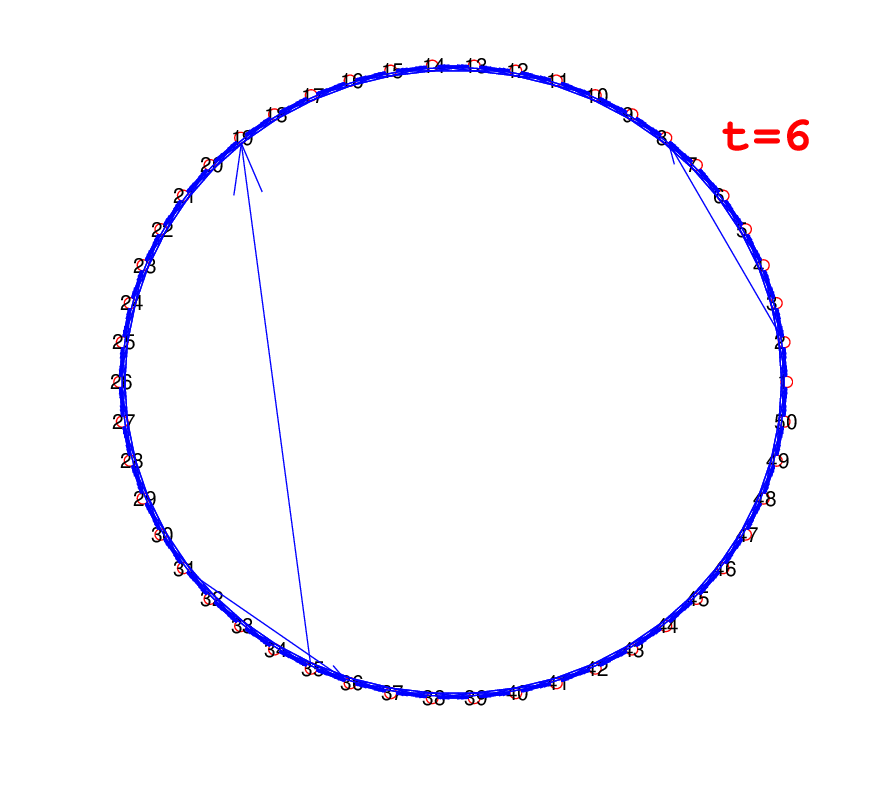}
\caption{The blinking model of shortcuts connections. Probability of switchings $p=0.04$, the switching time step $\tau=1$.}
\label{fig:2}
\end{figure}

In this section, a numerical example is given to demonstrate the effectiveness of the presented results on synchronization of LCODEs with time-varying couplings. The Lyapunov exponents are be computed numerically. By this way, we can verify the the synchronization
criterion and analyze synchronization numerically. We use the R\"{o}ssler system \cite{Sti,OER} as the node dynamics
\begin{eqnarray}\label{ross}
\left\{\begin{array}{lll}\dot{x}_{1}(t)=-{x}_{2}(t)-{x}_{3}(t)\\
\dot{x}_{2}(t)={x}_{1}(t)+a{x}_{2}(t)\\
\dot{x}_{3}(t)=b+{x}_{3}(t)({x}_{1}(t)-c)\end{array}\right.,
\end{eqnarray}
where $a=0.165, b=0.2$, and $c=10$. Figure \ref{fig:1} shows the dynamical behaviors of the R\"{o}ssler system (\ref{ross}) with random initial value in $[0,1]$ that includes a chaotic attractor \cite{Sti,OER}.

The network with time-varying topology we used here is NW small-world network with a time-varying coupling which was introduced as the blinking model  in \cite{Bel1,MEJ}. The time-varying network model algorithm is presented as follows. We divide the time axis into intervals of length $\tau$, in each interval: (a) Begin with a nearest neighbor coupled
network consisting of $m$ nodes arranged in a ring, where each node $i$ is adjacent to its $2k$-nearest neighbor nodes; (b) Add a connection between each pair of nodes with probability $p$, which usually is a random number between $[0,0.1]$;  For more details, we refer the readers to \cite{Bel1}. Figure \ref{fig:2} shows the time-varying structure of shortcut connections in the blinking model with $m=50$ and $k=3$.

In this example, the parameters
are set as $m=50$, $k = 3$, $\tau=1$, and $p=0.04$. Then small-world
network can be generated with the coupling {\it graph Laplacian} $\mathcal L (G(t))=-L(t)$. The dynamical network system can be describe as following:
\begin{eqnarray}\label{ms}
\left\{\begin{array}{lll}\dot{x}^{i}_{1}(t)=-{x}^{i}_{2}(t)-{x}^{i}_{3}(t)+\sigma\sum_{j=1}^{m}l_{ij}(t)x^{j}_{1}(t)\\
\dot{x}^{i}_{2}(t)={x}^{i}_{1}(t)+a{x}^{i}_{2}(t)+\sigma\sum_{j=1}^{m}l_{ij}(t)x^{j}_{2}(t)\\
\dot{x}^{i}_{3}(t)=b+{x}^{i}_{3}(t)({x}^{i}_{1}(t)-c)+\sigma\sum_{j=1}^{m}l_{ij}(t)x^{j}_{3}(t)\end{array}\right., ~i=\onetom.
\end{eqnarray}

\begin{figure}[hbt]
\centering
\includegraphics[width=1\textwidth,height=.6\textwidth]{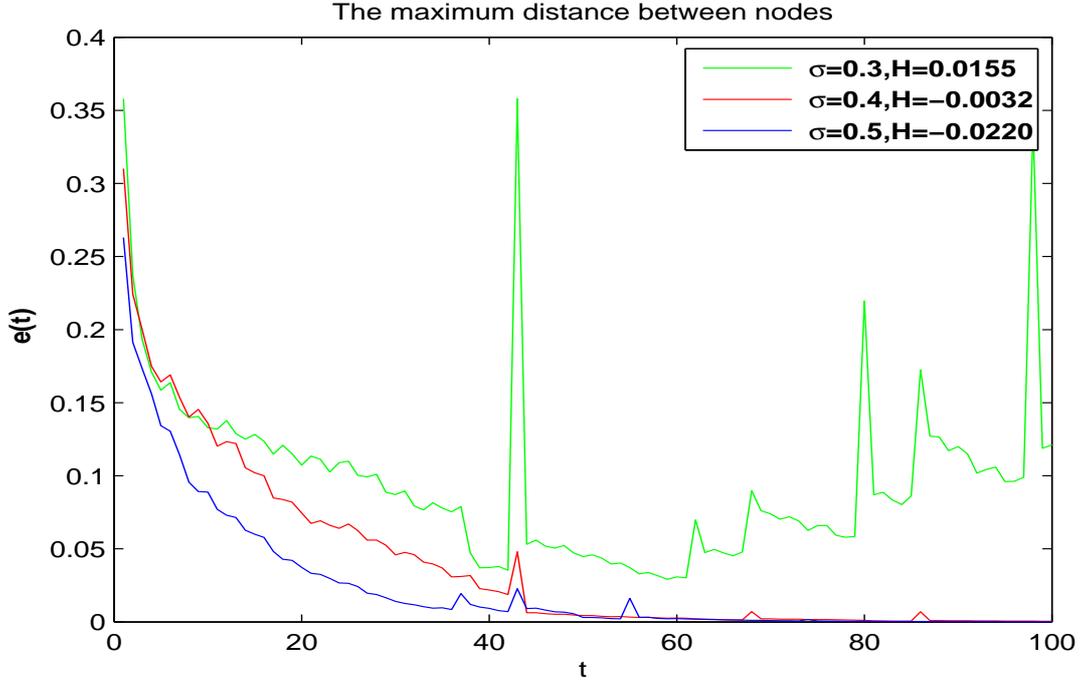}
\caption{Convergence of the maximum distance between nodes with different coupling strength $\sigma$.}
\label{fig:4}
\end{figure}

Let $e(t)=\max_{1\le i<j\le50}\|x^{i}(t)-x^{j}(t)\|$ denotes the maximum distance between nodes at time $t$. Let $E=\int_{T}^{T+R}e(t)dt$, for some sufficiently large $T>0$ and $R>0$. Let $H=\mu+\varsigma$ defined in Corollary 4. Figure \ref{fig:4} shows convergence of the maximum distance between nodes during the topology evolution with different coupling strength $\sigma$. It can be seen from Figure \ref{fig:4} that the dynamical network system (\ref{ms}) can be synchronized with $\sigma=0.4$ and $\sigma=0.5$.

\begin{figure}[hbt]
\centering
\includegraphics[width=1\textwidth,height=.6\textwidth]{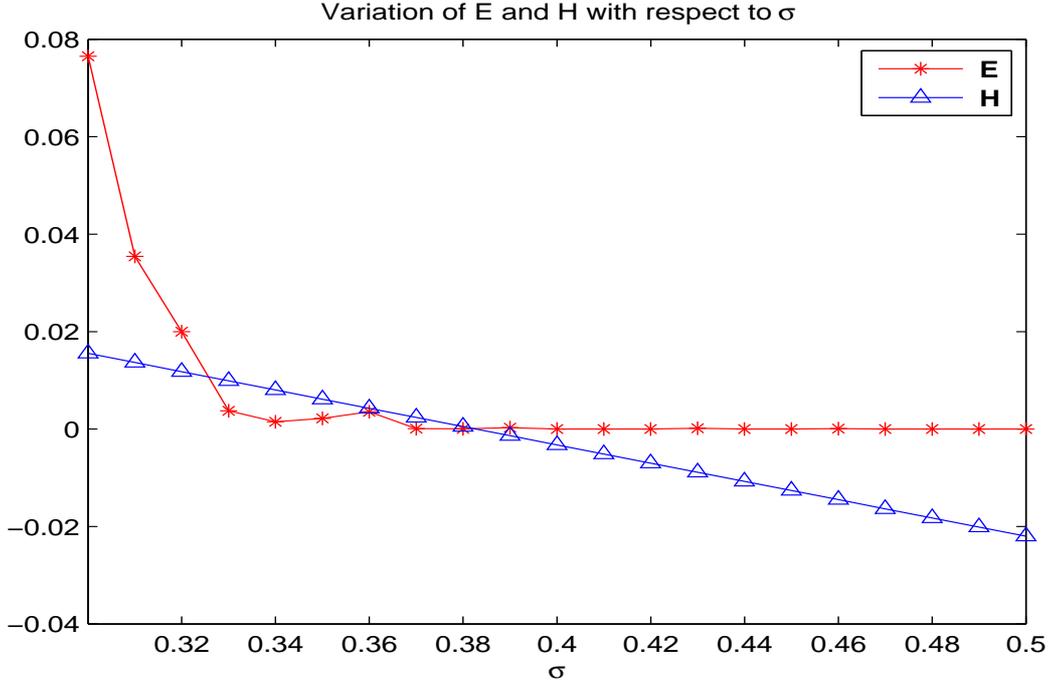}
\caption{Variation of $e$ and $H$ with respect to $\sigma$ for the blinking topology.}
\label{fig:5}
\end{figure}

We pick the evolution time length to be $200$. Let $T=190$, $R=10$. And choose initial state randomly from the interval $[0,1]$. Figure \ref{fig:5} shows the variation of $E$ and $H$ with respect to the coupling strength $\sigma$. It can be seen that the parameter (coupling strength $\sigma$) region where $H$ is negative coincides with
that of synchronization , i.e., where $E$ is near zero. This verified the theoretical result (Corollary 4). In addition, we find that $\sigma\approx 0.38$ is the threshold for synchronizing the coupled systems in this case.

\section{Conclusions}
In this paper, we present a theoretical framework for synchronization analysis of general coupled differential dynamical systems. The extended Hajnal diameter is introduced to measure the synchronization. The coupling between nodes is time-varying in both the network structure and the reaction dynamics. Inspired by the approaches in \cite{LAJ,LAJ1}, we show that the Hajnal diameter of the linear system induced by the time-varying coupling matrix and the largest Lyapunov exponent of the synchronized system play the key roles in synchronization analysis of LCODEs. These results extend synchronization analysis of discrete-time network in \cite{LAJ} to continuous-time case. As an application, we obtain a sufficient condition ensuring directed time-varying graph reaching consensus which is very general, and the way we get this result is different from \cite{Mor}. An example of numerical simulation is provided to show the effectiveness the theoretical results. Additional contributions of this paper are that we explicitly show that the largest projection Lyapunov exponent, the Hajanal diameter and the largest Lyapunov exponent of the transverse space equal in coupled differential systems (see lemma \ref{lem6} and lemma \ref{lem7}), which was proved in \cite{LAJ} for couple discrete-time systems.

\section*{Acknowledgements}
This work is jointly supported by the National Key Basic Research and Development Program (No. 2010CB731403), the National Natural Sciences Foundation of China under Grant (Nos. 61273211 and 61273309), the Shanghai Rising-Star Program (No. 11QA1400400), the
Marie Curie International Incoming Fellowship from the European Commission
under Grant FP7-PEOPLE-2011-IIF-302421, and the Laboratory of Mathematics for Nonlinear Science, Fudan University.

\section*{Appendices}
{\bf Proof of Lemma \ref{lem2}}: Let $U$ be a bounded open neighborhood of $A$
satisfying $\bigcap\limits_{t\ge 0}\vartheta^{(t)}\bar{U}=A$ and
$U_{t}=\{x\in R^{n}:\vartheta^{(\tau)}x\in U,0\le\tau\le t\}$. This
implies $U_{t}\supset U_{t'}$ if $t'\ge t\ge 0$, $U_{t}$ is an
open set due to the continuity of the semi-flow $\vartheta^{(t)}$,
and $\vartheta^{(\delta)} U_{t}\subset U_{t-\delta}$ for all
$t\ge \delta\ge 0$. Let $V=\bigcap\limits_{t\ge 0}U_{t}$. We claim
that there exists $t_{0}\ge 0$ such that $V=U_{t}$ for all
$t\ge t_{0}$.

For any $\delta>0$, let $t_{n}=n\delta$ and $U_{n}=U_{t_{n}}$. We
can conclude that $V=\bigcap\limits_{n=1}^{\infty}U_{n}$. We will
prove in the following that there exists $n_{0}$ such that
$V=U_{n_{0}}$. Otherwise, there always exists $x_{n}\in
U_{n}\setminus U_{n+1}$ for $n\ge 0$. Let
$y_{n}=\vartheta^{(t_{n+1})}x_{n}$. We have (i)
$y_{n}\in\bigcap\limits_{k=0}^{n}\vartheta^{(t_{k})}\bar{U}$ and
(ii) $y_{n}\notin U$. For any limit point $\hat{y}$ of $y_{n}$,
$\hat{y}$ can be either finite or infinite. For both cases,
$\hat{y}\notin U$ which implies $\hat{y}\notin A$. However, the
claim (i) implies that $\hat{y}\in A$, which contradicts with the
claim (ii). This completes the proof by letting $W=\bar{V}$.

{\bf Proof of Lemma \ref{lem5}}: (a) For any initial condition
with the form $\delta x_{0}=\mathbf 1_{m}\otimes u_{0} $, the
solution of (\ref{cauchy}) can be $U(t,t_{0},s_{0})\bigg(\mathbf
1_{m}\otimes u_{0}\bigg)=\mathbf 1_{m}\otimes
\breve{U}(t,t_{0},s_{0})u_{0}$ according to Lemma \ref{lem1}. This
implies the first claim in this lemma.

(b) According to Lemma \ref{lem2}, there exists $K_{1}>0$ such
that $s(t)$, the solution of (\ref{syn}),  satisfies $\|s(t)\|\le K_{1}$ for all
$s_{0}\in W$ and $t\ge 0$. So, there exists $K>0$ such that
$\|DF^{t}(s(t))\|\le K$ according to the $3$-th item of the
assumption $\mathbf A1$. Write the solution of (\ref{cauchy})
$\delta x(t)=U(t,t_{0},s_{0})\delta x_{0}$ as
\begin{eqnarray*}
\delta x(t+t_{0})=\delta
x_{0}+\int_{t_{0}}^{t+t_{0}}DF^{\tau}(s(\tau))\delta x(\tau)d\tau.
\end{eqnarray*}
Then,
\begin{eqnarray*}
\|\delta x(t+t_{0})\|\le\|\delta
x_{0}\|+\int_{t_{0}}^{t+t_{0}}\|DF^{\tau}(s(\tau))\|\|\delta
x(\tau)\|d\tau\le\|\delta
x_{0}\|+K\int_{0}^{t}\|\delta
x(\tau+t_{0})\|d\tau.
\end{eqnarray*}
According to Lemma \ref{Grounwell}, we have $\|\delta x(t+t_{0})\|\le \|\delta
x_{0}\|+K\int_{0}^{t}\|\delta
x_{0}\|e^{(t-\tau)K}d\tau=
e^{Kt}\|\delta x_{0}\|$. This implies that
$\|U(t+t_{0},t_{0},s_{0})\|\le e^{Kt}$ for all $s_{0}\in W$
and $t_{0}\ge 0$.

For any $s_{0},s'_{0}\in W$, let $s(t)$ and $s'(t)$ be the
solution of the synchronized state equation (\ref{syn}) with
initial condition $s(t_{0})=s_{0}$ and $s'(t_{0})=s'_{0}$
respectively. We have
\begin{align*}
s(t+t_{0})-s'(t+t_{0})=\int_{t_{0}}^{t+t_{0}}[f(s(\tau))-f(s'(\tau))]d\tau+s(t_{0})-s'(t_{0}),\\
\|s(t+t_{0})-s'(t+t_{0})\|\le\|s(t_{0})-s'(t_{0})\|+K\int_{t_{0}}^{t+t_{0}}\|s(\tau)-s'(\tau)\|d\tau.\\
\end{align*}
By Lemma \ref{Grounwell}, we have $\|s(t+t_{0})-s'(t+t_{0})\|\le e^{K
t}\|s_{0}-s'_{0}\|$ for all $t_{0},t\ge 0$ and $s_{0},s'_{0}\in
W$. Also, according to the $4$-th item of assumption $\mathbf
A1$, there must exist $K_{2}>0$ such that
$\|DF^{t}(s(t))-DF^{t}(s'(t))\|\le K_{2}\|s(t)-s'(t)\|$ for
all $t\ge 0$ and $s_{0},s'_{0}\in W$. Then, let $\delta
x(t)=U(t,t_{0},s_{0})\delta x_{0}$, $\delta
y(t)=U(t,t_{0},s_{0}')\delta x_{0}$, and $v(t)=\delta x(t)-\delta
y(t)$. We have
\begin{align*}
v(t+t_{0})&=\int_{t_{0}}^{t+t_{0}}\bigg[DF^{\tau}(s(\tau))\delta
x(\tau)-DF^{\tau}(s'(\tau))\delta y(\tau)\bigg]d\tau\\
&=\int_{t_{0}}^{t+t_{0}}[DF^{\tau}(s(\tau))-DF^{\tau}(s'(\tau))]\delta
x(\tau)d\tau+\int_{t_{0}}^{t+t_{0}}DF^{\tau}(s'(\tau))v(\tau)d\tau,\\
\|v(t+t_{0})\|&\le\int_{t_{0}}^{t+t_{0}}\bigg[\|DF^{\tau}(s(\tau))-DF^{\tau}(s'(\tau))\|\|\delta
x(\tau)\|+\|DF^{\tau}(s'(\tau))\|\|v(\tau)]\|\bigg]d\tau\\
&\le K_{2}\int_{0}^{t}e^{2K\tau}d\tau\|\delta
x_{0}\|\|s_{0}-s'_{0}\|+K\int_{t_{0}}^{t_{0}+t}\|v(\tau)\|d\tau.\\
\end{align*}
According to Lemma \ref{Grounwell},
\begin{eqnarray*}
\|v(t+t_{0})\|\le\bigg[\frac{K_{2}(e^{2Kt}-e^{Kt})}{K}\bigg]\|\delta
x_{0}\|\|s_{0}-s'_{0}\|
\end{eqnarray*}
This implies
\begin{eqnarray*}
\|U(t+t_{0},t_{0},s_{0})-U(t+t_{0},t_{0},s'_{0})\|\le\bigg[\frac{K_{2}(e^{2Kt}-e^{Kt})}{K}\bigg]\|s_{0}-s'_{0}\|
\end{eqnarray*}
for all $s_{0},s'_{0}\in W$. This completes the proof.

{\bf Proof of Lemma \ref{lem6}}: We define the projection joint
spectral radius as follows:
\begin{eqnarray*}
\rho_{P}(D\mathcal
F,s_{0})=\overline{\lim\limits_{t\rightarrow\infty}}\sup\limits_{t_{0}\ge
0}\|\tilde{U}(t,t_{0},s_{0})\|^{\frac{1}{t}}
\end{eqnarray*}
First, we will prove that ${\rm diam}(D\mathcal
F,s_{0})=\rho_{P}(D\mathcal F,s_{0})$. For any
$d>\rho_{P}(D\mathcal F,s_{0})$, there exists $T\ge 0$ such that
$\|\tilde{U}(t+t_{0},t_{0},s_{0})\|\le d^{t}$ for all
$t_{0}\ge 0$ and $t\ge T$. This implies that
\begin{align*}
&\bigg\|P^{-1}U(t+t_{0},t_{0},s_{0})P-\left[\begin{array}{l}I_{n}\\0\\\vdots\\0\end{array}\right]
\bigg[P_{0}^{\top}\breve{U}(t+t_{0},t_{0},s_{0})P_{0},\alpha(t+t_{0},t_{0},s_{0})\bigg]\bigg\|\\
&=\bigg\|\left[\begin{array}{ll}0&0\\
0&\tilde{U}(t+t_{0},t_{0},s_{0})\end{array}\right]\bigg\|
\le C_{1}d^{t}
\end{align*}
for some $C_{1}>0$, all $t_{0}\ge 0$, and all $t\ge T$.
Thus, there exist some $C_{2}>0$ and some matrix function
$q(t)\in\mathbb R^{n,nm}$ such that
\begin{align*}
&\|U(t+t_{0},t_{0},s_{0})-\mathbf 1_{m}\otimes q(t)\|\\
&=\bigg\|U(t+t_{0},t_{0},s_{0})-P\left[\begin{array}{l}I_{n}\\0\\\vdots\\0\end{array}\right]
\bigg[P_{0}^{\top}\breve{U}(t+t_{0},t_{0},s_{0})P_{0},\alpha(t+t_{0},t_{0},s_{0})\bigg]P^{-1}v\bigg\|\le C_{2}d^{t}
\end{align*}
for all $t_{0}\ge 0$ and $t\ge T$, where $q(t)\in \mathbb R^{n,nm}$ denotes a appropriate matrix we omit its accurate expression. So, we can conclude that
${\rm diam}(U(t+t_{0},t_{0},s_{0}))\le C_{3}d^{t}$ for some
$C_{3}>0$, all $t_{0}\ge 0$, and $t\ge T$. This implies that ${\rm
diam}(D\mathcal F,s_{0})\le d$, i.e., ${\rm diam}(D\mathcal
F,s_{0})\le\rho_{P}(D\mathcal F,s_{0})$ due to the arbitrariness
of $d\ge\rho_{P}(D\mathcal F, s_{0})$. Conversely, for any $d>{\rm
diam}(D\mathcal F,s_{0})$, there exists $T>0$ such that
\begin{eqnarray*}
\|U(t+t_{0},t_{0},s_{0})-\mathbf 1_{m}\otimes U_{1}\|\le
C_{4}d^{t}
\end{eqnarray*}
 for some $C_{4}>0$, all $t_{0}\ge 0$, and $t\ge T$, where
$U_{1}=[U_{11},U_{12},\cdots,U_{1m}]$ the first $n$ rows of $U(t+t_{0},t_{0},s_{0})$. Then,
\begin{align*}
&\|P^{-1}U(t+t_{0},t_{0},s_{0})P-P^{-1}\mathbf 1_{m}\otimes U_{1}P\|=\bigg\|P^{-1}U(t+t_{0},t_{0},s_{0})P-\left[\begin{array}{ll}\gamma(t)&\beta(t)\\
0&0\end{array}\right]\bigg\|\\
&=\bigg\|\left[\begin{array}{ll}0&\beta(t)\\
0&\tilde{U}(t+t_{0},t_{0},s_{0})\end{array}\right]\bigg\|\le C_{5}d^{t}
\end{align*}
for some $C_{5}>0$, all $t_{0}\ge 0$, and $t\ge T$, where
$\gamma(t)=P_{0}^{\top}\breve{U}(t,t_{0},s_{0})P_{0}\in\mathbb R^{n,n}$ and $\beta(t)\in\mathbb
R^{n,n(m-1)}$ denotes a appropriate matrix we omit its accurate expression. This implies that
$\|\tilde{U}(t+t_{0},t_{0},s_{0})\|\le C_{6}d^{t}$ holds for some
$C_{6}>0$, all $t_{0}\ge 0$, and $t\ge T$. Therefore, we can
conclude that $\rho_{P}(D\mathcal F,s_{0})\le d$. So,
$\rho_{P}(D\mathcal F,s_{0})={\rm diam}(D\mathcal F,s_{0})$.

Second, it is clear that $\log\rho_{P}(D\mathcal
F,s_{0})\ge\lambda_{P}(D\mathcal F,s_{0})$. We will prove that
$\log\rho_{P}(D\mathcal F,s_{0})=\lambda_{P}(D\mathcal
F,s_{0})$. Otherwise, there exists some $r,r_{0}>0$ satisfying
$\rho_{P}(D\mathcal F,s_{0})>r>r_{0}>e^{\lambda_{P}(D\mathcal
F,s_{0})}$. If so, there exist an sequence $t_{k}\uparrow\infty$
as $k\rightarrow\infty$, $t_{0}^{k}\ge 0$, and $v_{k}\in\mathbb
R^{n(m-1)}$ with $\|v_{k}\|=1$ such that
$\|\tilde{U}(t_{k}+t_{0}^{k},t_{0}^{k},s_{0})v_{k}\|>r^{t_{k}}$
for all $k\in\mathcal N$. Then, there exists a subsequence
$v_{k_{l}}$ with
$\lim\limits_{l\rightarrow\infty}v_{k_{l}}=v^{*}$. Let
$\{e_{1},e_{2},\cdots,e_{n(m-1)}\}$ be a normalized orthogonal
basis of $\mathbb R^{n(m-1)}$. And, let
$v_{k_{l}}-v^{*}=\sum\limits_{j=1}^{n(m-1)}\xi_{j}^{k_{l}}e_{j}$.
We have $\lim\limits_{l\rightarrow\infty}\xi_{j}^{k_{l}}=0$
for all $j=1,\cdots,n(m-1)$. Thus, there exists $L>0$ such that
\begin{align*}
\|\tilde{U}(t_{k_{l}}+t_{0}^{k_{l}},t_{0}^{k_{l}},s_{0})v^{*}\|&\ge\|\tilde{U}(t_{k_{l}}+t_{0}^{k_{l}},t_{0}^{k_{l}},s_{0})v_{k_{l}}\|-\|\tilde{U}(t_{k_{l}}+t_{0}^{k_{l}},t_{0}^{k_{l}},s_{0})(v_{k_{l}}-v^{*})\|\\
&\ge r^{t_{k_{l}}}-\sum\limits_{j=1}^{n(m-1)}|\xi_{j}^{k_{l}}|\|\tilde{U}(t_{k_{l}}+t_{0}^{k_{l}},t_{0}^{k_{l}},s_{0})e_{j}\|\ge r^{t_{k_{l}}}-r_{0}^{t_{k_{l}}}>r_{0}^{t_{k_{l}}}\\
\end{align*}
for all $l\ge L$. This implies $e^{\lambda(D_{P}\mathcal
F,v^{*},s_{0})}\ge r_{0}$ which contradicts with
$e^{\lambda_{P}(D\mathcal F,s_{0})}<r_{0}$. This implies
$\rho_{P}(D\mathcal F,s_{0})=e^{\lambda_{P}(D\mathcal F,s_{0})}$.
Therefore, we can conclude $\log{\rm diam}(D\mathcal
F,s_{0})=\lambda_{P}(\mathcal F,s_{0})$. The proof is completed.

{\bf Proof of Lemma \ref{lem7}}: Let $\tilde{\phi}=P^{-1}\phi$. We
have
\begin{align*}
\dot{\tilde{\phi}}=P^{-1}DF(s(t),\varrho^{(t)}\omega_{0})P\tilde{\phi}=\left[\begin{array}{ll}
P_{0}^{\top}\frac{\partial f}{\partial s}(s(t))P_{0}&\alpha(t)\\
0&\tilde{D}F(s(t),\varrho^{(t)}\omega_{0})\end{array}\right]\tilde{\phi}.
\end{align*}
Write
$\tilde{\phi}=\left[\begin{array}{l}y(t)\\z(t)\end{array}\right]$,
where $y(t)\in\mathbb R^{n}$. Then, we have
\begin{align*}
\left\{\begin{array}{l}\dot{z}(t)=\tilde{D}F(s(t),\varrho^{(t)}\omega_{0})z(t)\\
\dot{y}(t)=P_{0}^{\top}\frac{\partial f}{\partial
s}(s(t))P_{0}y(t)+\alpha(t)z(t)\end{array}\right..
\end{align*}
Thus, we can write its solution by
\begin{eqnarray*}
\left\{\begin{array}{l}z(t)=\tilde{U}(t)z_{0}\\
y(t)=P_{0}^{\top}\breve{U}(t)P_{0}y_{0}+\int_{0}^{t}P_{0}^{\top}\breve{U}(t)\breve{U}^{-1}(\tau)P_{0}
\alpha(\tau)\tilde{U}(\tau)z_{0}d\tau\end{array}\right..
\end{eqnarray*}

We write $\lambda_{P}(D\mathcal F, s_{0},\omega_{0})$,
$\lambda_{S}(D\mathcal F, s_{0},\omega_{0})$, and
$\lambda_{T}(D\mathcal F, s_{0},\omega_{0})$ by $\lambda_{P}$,
$\lambda_{S}$, and $\lambda_{T}$ respectively for simplicity.

{\bf Case 1}: $\lambda_{P}>\lambda_{S}$. We can conclude that
$\chi[z(t)]\le \lambda_{P}$ and
\begin{align*}
\chi[y(t)]\le\max\bigg\{\chi\bigg[P_{0}^{\top}\breve{U}(t)P_{0}y_{0}\bigg],~\chi\bigg[\int_{0}^{t}P_{0}^{\top}\breve{U}(t)\breve{U}^{-1}(\tau)P_{0}\alpha(\tau)\tilde{U}(\tau)z(0)d\tau\bigg]\bigg\}.
\end{align*}
From Cauchy-Buniakowski-Schwarz inequality, we have
\begin{align*}
&\chi\bigg[\bigg\|\int_{0}^{t}P_{0}^{\top}\breve{U}(t)\breve{U}^{-1}(\tau)P_{0}\alpha(\tau)\tilde{U}(\tau)d\tau\bigg\|\bigg]\\
&\le\chi\bigg[\bigg\{\int_{0}^{t}\|\breve{U}(t)\breve{U}^{-1}(\tau)\|^{2}d\tau\bigg\}^{1/2}\bigg]
+\chi\bigg[\bigg\{\int_{0}^{t}\|\alpha(\tau)\tilde{U}(\tau)\|d\tau\bigg\}^{1/2}\bigg].
\end{align*}
{\bf Claim 1}:
$\chi(\int_{0}^{t}\|\breve{U}(t)\breve{U}^{-1}(\tau)\|^{2}d\tau)\le
0$.

Considering the linear system
\begin{eqnarray}
\dot{u}(t)=\frac{\partial f}{\partial s}(s(t))u(t),\label{varsyn}
\end{eqnarray}
due to its regularity and the boundedness of its coefficients, there exists a Lyapunov transform $L(t)$
such that letting $u(t)=L(t)v(t)$, consider the transformed linear
system
\begin{eqnarray}
\dot{v}(t)=\bigg[L^{-1}(t)\frac{\partial f}{\partial
s}(s(t))L(t)-L^{-1}(t)\dot{L}(t)\bigg]v(t)=\breve{A}(t)v(t).
\end{eqnarray}
Let solution matrix
$\breve{V}(t)=(\breve{v}_{ij}(t))_{i,j=1}^{n}$,
$\breve{A}(t)=(\breve{a}_{ij}(t))_{i,j=1}^{n}$ which satisfies
that $\breve{A}(t)$ and $\breve{V}(t)$ are lower-triangular. And,
its Lyapunov exponents can be written as follows:
\begin{eqnarray*}
\sigma_{i}=\lim\limits_{t\rightarrow\infty}\frac{1}{t}\int_{0}^{t}\breve{a}_{ii}(\tau)d\tau,
\end{eqnarray*}
which are just the Lyapunov exponents of the regular linear system
(\ref{varsyn}), $i=\oneton$. We have
$\chi[\breve{v}_{ii}(t)]=\sigma_{i}$ and
\begin{eqnarray*}
\breve{v}_{k+1,k}(t)=e^{\int_{0}^{t}\breve{a}_{k+1,k+1}(\tau)d\tau}\int_{0}^{t}e^{-\int_{0}^{\tau}\breve{a}_{k+1,k+1}(\vartheta)d\vartheta}
\breve{a}_{k+1,k}(\tau)\breve{v}_{k,k}(\tau)d\tau.
\end{eqnarray*}
This implies
\begin{eqnarray*}
\chi[\breve{v}_{k+1,k}(t)]\le
\sigma_{k+1}-\sigma_{k+1}+0+\sigma_{k}=\sigma_{k}.
\end{eqnarray*}
By induction, we can conclude that $\chi[\breve{v}_{jk}(t)]\le
\sigma_{k}$ for all $j> k$. For $j<k$,
$\chi[\breve{v}_{jk}(t)]=-\infty$ due to the lower-triangularity
of the matrix $\breve{V}(t)$.

Considering the lower-triangular matrix
$\breve{V}^{-1}(t)=(\breve{w}_{ij})_{i,j=1}^{n}$, its transpose
$(\breve{V}^{-1}(t))^{\top}$ can be regarded as the solution
matrix of the adjoint system of (\ref{varsyn}):
\begin{eqnarray*}
\dot{w}(t)=-\breve{A}^{\top}(t)w(t),
\end{eqnarray*}
which is also regular. By the same arguments, we can conclude that
$\chi[\breve{w}_{kk}]=-\sigma_{k}$ for all $k=\oneton$,
$\chi[\breve{w}_{jk}]\le-\sigma_{k}$ for all $k>j$, and
$\chi[\breve{w}_{jk}]=-\infty$ for all $k<j$. Therefore, for
each $i>j$,
\begin{align*}
\max_{i,j}\chi\bigg[\int_{0}^{t}|\breve{U}(t)\breve{U}^{-1}(\tau)|_{ij}d\tau\bigg]&\le\max_{i,j}\chi\bigg[\int_{0}^{t}|\breve{V}(t)\breve{V}^{-1}(\tau)|_{ij}d\tau\bigg]\le\max_{i,j}\chi\bigg[\int_{0}^{t}\sum\limits_{j\le
k\le
i}|\breve{v}_{ik}(t)\breve{w}_{kj}(\tau)|d\tau\bigg]\\
&\le\max_{i,j}\max\limits_{j\le k\le
i}\chi\bigg[\int_{0}^{t}|\breve{v}_{ik}(t)\breve{w}_{kj}(\tau)|d\tau\bigg]\le\max_{i,j}\max\limits_{j\le
k\le i}(\sigma_{k}-\sigma_{k})=0.
\end{align*}
This implies that
$\chi[\int_{0}^{t}\|\breve{U}(t)\breve{U}^{-1}(\tau)d\tau\|^{2}d\tau]\le
0$.

Noting that
\begin{eqnarray*}
\chi\bigg[\int_{0}^{t}\|\alpha(\tau)\tilde{U}(\tau)\|_{2}^{2}d\tau\bigg]\le
 \chi\bigg[\|\alpha(t)\tilde{U}(t)\|^{2}\bigg]\le
2\lambda_{P}.
\end{eqnarray*}
So, $\chi[y(t)]\le\max\{\lambda_{S},\lambda_{P}\}=\lambda_{P}$.
This leads to $\chi[\tilde{\phi}(t)]\le \lambda_{P}$. This implies
that $\lambda_{P}=\max\{\lambda_{S},\lambda_{T}\}$. Thus,
$\lambda_{P}=\lambda_{T}$ can be concluded due to
$\lambda_{P}>\lambda_{S}$.

{\bf Case 2}: $\lambda_{P}<\lambda_{S}$. For any $\epsilon$ with
$0<\epsilon<\frac{\lambda_{S}-\lambda_{P}}{3}$, there exists $T>0$
such that
\begin{eqnarray*}
\|\breve{U}^{-1}(\tau)\|\le e^{(-\lambda_{S}+\epsilon)\tau},~\|\alpha(\tau)\|\le e^{\epsilon\tau},~\|\tilde{U}(\tau)\|\le e^{(\lambda_{P}+\epsilon)\tau}
\end{eqnarray*}
for all $t\ge T$. Define the subspace of $\mathbb R^{nm}$:
\begin{eqnarray*}
V=\bigg\{\left[\begin{array}{ll}y\\z\end{array}\right]:~y=-\int_{0}^{\infty}P_{0}^{\top}\breve{U}^{-1}(\tau)P_{0}\alpha(\tau)
\tilde{U}(\tau)d\tau z\bigg\},
\end{eqnarray*}
which is well defined due to
$\|P_{0}^{\top}\breve{U}^{-1}(\tau)P_{0}\alpha(\tau)\tilde{U}(\tau)\|\le
e^{(3\epsilon-\lambda_{S}+\lambda_{P})\tau}\in L([T,+\infty))$.
For each $\tilde{\phi}(t)$ with initial condition
$\left[\begin{array}{ll}y\\z\end{array}\right]\in V$, we have
$\chi[z(t)]\le\lambda_{P}$ and
\begin{align*}
&\chi[y(t)]\\&=\chi\bigg[\bigg\{-P_{0}^{-1}\breve{U}(t)P_{0}\int_{0}^{\infty}P_{0}^{\top}\breve{U}^{-1}(\tau)P_{0}\alpha(\tau)\tilde{U}(\tau)d\tau
+P_{0}^{\top}\breve{U}(t)P_{0}\int_{0}^{t}P_{0}^{\top}\breve{U}^{-1}(\tau)P_{0}\alpha(\tau)\tilde{U}(\tau)d\tau\bigg\}z\bigg]\\
&=\chi\bigg[-P_{0}^{\top}\breve{U}(t)\int_{t}^{\infty}\breve{U}^{-1}(\tau)P_{0}\alpha(\tau)\tilde{U}(\tau)d\tau
z\bigg]\le\lambda_{P}
\end{align*}
according to the arguments above. Thus, we have $\max\limits_{u\in
V}\lambda(D\mathcal F,u,s_{0},\omega_{0})=\lambda_{P}$. Since
$\dim (V)=n(m-1)$, $V$ define the transverse space and
$\lambda_{T}=\lambda_{P}$. This completes the proof.

{\bf Proof of Lemma \ref{lem8}}: Since $L(t)$ satisfies assumption
$\mathbf A_{4}$, if the initial condition is $u(t_{0})={\bf
1}_{m}$, then the solution must be $u(t)={\bf 1}_{m}$, which
implies that each row sum of $V(t,t_{0})$ is one. Then, we will
prove all elements in $V(t,t_{0})$ are nonnegative. Consider the
$ith$ column of $V(t,t_{0})$ denoted by $V^{i}(t,t_{0})$ which can
be regarded as the solution of the following equation:
\begin{eqnarray*}
\left\{\begin{array}{l}\dot{u}=\sigma L(t)u\\
u(t_{0})=e_{i}^{m}\end{array}\right.
\end{eqnarray*}
For any $t\ge t_{0}$, if $i_{0}=i_{0}(t)$ is the index with
$u_{i_{0}}(t)=\min\limits_{i=\onetom}u_{i}(t)$, we have
$\dot{u}_{i_{0}}(t)=\sum\limits_{j=1}^{m}\sigma l_{i_{0}j}(u_{j}(t)-u_{i_{0}}(t))\ge
0$. This implies that $\min\limits_{i=\onetom}u_{i}(t)$ is always
nondecreasing for all $t\ge t_{0}$. Therefore, $u_{i}(t)\ge 0$
holds for all $i=\onetom$ and $t\ge t_{0}$. We can conclude that
$V(t,t_{0})$ is a stochastic matrix. The proof is completed.

{\bf Proof of Lemma \ref{lem9}}: Consider the following Cauchy
problem:
\begin{eqnarray*}
\left\{\begin{array}{l}
\dot{u}_{i}(t)=\sum\limits_{j=1}^{m}\sigma l_{ij}(t)u_{j}(t)\\
u_{i}(t_{0})=\left\{\begin{array}{ll}1&i=k\\
0&{\rm otherwise}\end{array}\right.\end{array}\right.,~i=\onetom.
\end{eqnarray*}
Noting that $\dot{u}_{k}(t)\ge \sigma l_{kk}u_{k}$, we have $u_{k}(t)\ge
e^{-M_{1}(t-t_{0})}$. For each $i\ne k$, since $u_{i}(t)\ge 0$
for all $i=\onetom$ and $t\ge t_{0}$,  we have
\begin{align*}
u_{i}(t)&=\sum\limits_{j\ne
i}\int_{t_{0}}^{t}e^{\int_{\tau}^{t}\sigma l_{ii}(\vartheta)d\vartheta}\sigma l_{ij}(\tau)u_{j}(\tau)d\tau\ge\int_{t_{0}}^{t}e^{\int_{\tau}^{t}\sigma l_{ii}(\vartheta)d\vartheta}\sigma l_{ik}(\tau)u_{k}(\tau)d\tau\\
&\ge\int_{i_{0}}^{t}e^{-M_{1}(t-\tau)}e^{-M_{1}(\tau-t_{0})}\sigma l_{ik}(\tau)d\tau=e^{-M_{1}(t-t_{0})}\int_{t_{0}}^{t}\sigma l_{ik}(\tau)d\tau.
\end{align*}
So, if there exists a $\delta$-edge from vertex $j$ to $i$ across
$[t_{0},t_{0}+T]$, then we have $v_{ij}(t_{0}+T,t_{0})\ge
e^{-M_{1}T}\delta$. Let
$\delta_{2}=\min\{e^{-M_{1}T},e^{-M_{1}T}\delta\}$. We can see
that $V(t,t_{0})$ has a $\delta_{2}$ spanning tree across any
$T$-length time interval. Therefore, according to \cite{Wolf, Ose},
there exist $\delta_{1}>0$ and $T_{1}=(m-1)T$ such that
$V(t,t_{0})$ is $\delta_{1}$ scrambling across any $T_{1}$-length
time interval. The Lemma is proved.

\end{document}